\documentclass[letterpaper]{article}
\usepackage[utf8]{inputenc}
\usepackage{amsmath}
\usepackage{amssymb,amsfonts,textcomp}
\usepackage[T1]{fontenc}
\usepackage[english,italian]{babel}
\usepackage{color}
\usepackage{array}
\usepackage{hhline}
\usepackage{hyperref}
\usepackage[pdftex]{graphicx}

\newcommand\textstyleFootnoteSymbol[1]{\textsuperscript{#1}}
\newcommand\textstyleRimandonotaapiiii[1]{\textsuperscript{#1}}

\makeatletter
\newcommand\ps@Standard{
  \renewcommand\@oddhead{}
  \renewcommand\@evenhead{}
  \renewcommand\@evenfoot{\@oddfoot}
  \renewcommand\thepage{\arabic{page}}
}
\makeatother
\newcommand\liststyleWWNumv{%
\renewcommand\labelitemi{{}-}
\renewcommand\labelitemii{o}
\renewcommand\labelitemiii{[F0A7?]}
\renewcommand\labelitemiv{[F0B7?]}
}
\newcommand\liststyleWWNumi{%
\renewcommand\labelitemi{{}-}
\renewcommand\labelitemii{o}
\renewcommand\labelitemiii{[F0A7?]}
\renewcommand\labelitemiv{[F0B7?]}
}
\title{}
\author{Struppa, Daniele}
\date{2018-08-07}
\begin{document}
\clearpage\setcounter{page}{1}\pagestyle{Standard}
\section[The Mathematics of Painting: the Birth of Projective Geometry in the Italian Renaissance ]{\textbf{\center The
Mathematics of Painting: \\
the Birth of Projective Geometry in the Italian Renaissance }}
\section[]{\bfseries }

\begin{center}
{\large Graziano Gentili, Luisa Simonutti and Daniele C. Struppa}\footnotetext{The first author was partially supported by INdAM and by Chapman University.
The second author was partially supported by ISPF-CNR and by Chapman University}
\end{center}

\bigskip

{\raggedleft\selectlanguage{english}
«Porticus aequali quamvis est denique ductu
\par}

{\raggedleft\selectlanguage{english}
\foreignlanguage{italian}{stansque in perpetuum paribus suffulta columnis,}
\par}

{\raggedleft\selectlanguage{english}
\foreignlanguage{italian}{longa tamen parte ab summa cum tota videtur,}
\par}

{\raggedleft\selectlanguage{english}
\foreignlanguage{italian}{paulatim trahit angusti fastigia coni,}
\par}

{\raggedleft\selectlanguage{english}
\foreignlanguage{italian}{tecta solo iungens atque omnia dextera laevis}
\par}

{\raggedleft\selectlanguage{english}
donec in obscurum coni conduxit acumen.»
\par}

\bigskip

{\raggedleft\selectlanguage{english}
Titus Lucretius Carus, \textit{De rerum natura}, IV 426-431
\par}

\bigskip

{\selectlanguage{english}
\textbf{Abstract.} We show how the birth of perspective painting in the Italian Renaissance led to a new way of
interpreting space that resulted in the creation of projective geometry. Unlike other works on this subject, we
explicitly show how the craft of the painters implied the introduction of new points and lines (points and lines at
infinity) and their projective coordinates to complete the Euclidean space to what is now called projective space. We
demonstrate this idea by looking at original paintings from the Renaissance, and by carrying out the explicit analytic
calculations that underpin those masterpieces.}

\bigskip

{\selectlanguage{english}
\textbf{Keywords.} Renaissance, Piero della Francesca, painting, perspective, analytic projective geometry, points and
lines at infinity.}

\section{1. Introduction}
{\selectlanguage{english}
The birth of projective geometry through the contribution of Italian Renaissance painters is a topic that has originated
a large and very interesting bibliography, some of which is referred to in this article. Most of the existing
literature dwells on the evolution of the understanding of the techniques that painters and artists such as Leon
Battista Alberti and Piero della Francesca developed to assist them (and other painters) in creating realistic
representations of scenes. These techniques, of course, are a concrete translation of ideas that slowly germinated and
were only later completely developed into a new branch of geometry that goes under the name of projective geometry. }

{\selectlanguage{english}
The point of view that we are taking in this article, however, is to strengthen the linkage between the pictorial ideas
and the mathematical underpinnings. More to the point, the entire architecture of prospective painting consists in
realizing that the space of vision cannot be represented through the usual Euclidean space, but requires the inclusion
of new geometrical objects that, properly speaking, do not exist in the Euclidean space. We are referring here to what
mathematicians call improper points and improper lines or, with a more suggestive term, points and lines at infinity.
Unlike most other studies, for example [4] and [12], we use here the approach and the terminology from analytic
projective geometry, rather than the proportion theory from Euclidean geometry, by introducing the notion of projective
coordinates. Just as the birth of projective geometry was stimulated by pictorial necessities, we show here how the
language of this new geometry can be applied to those necessities.}

{\selectlanguage{english}
There are two main reasons for this approach. On one hand, we believe the projective terminology allows a simpler way to
treat the technical task at hand, namely the identification of the technical processes that a painter needs to
represent a scene. But there is a deeper reason: perspective is not simply a technique; rather it is a radical change
of perspective (pun intended) on what space is. In order to formally perfect the process of representation, the
mathematicians had to introduce new objects, new points, new lines, new planes. It is by introducing these objects that
mathematicians were able to create a logically consistent view of the pictorial space, that allowed them a formally
unimpeachable process through which what we see can be translated into what we draw. The new line, plane, space (which
are now the projective line, the projective plane, the projective space) resemble (and contain) the old Euclidean line,
plane, space, but perfect the nature of their properties. So, for example, while in the Euclidean plane we say that any
two distinct lines intersect in a point \textit{unless} they are parallel, in the new projective plane we can say that
any two distinct lines intersect in a point, without exception. Projective geometry is not just a new and useful
technique, it is a radically different way of representing the space around us.}

{\selectlanguage{english}
We should add a couple of notes for the reader. Projective geometry is born of the necessity to understand the
phenomenon of apparent intersection between parallel lines, and most of our article is devoted to this aspect. However,
once the mathematics is clear, projective geometry allows the study of much more complex situations. For example, the
same techniques that we illustrate in our article, can be utilized to determine how to represent the halo of a saint,
or the shadow of a lamp against the wall of a church. This topic goes beyond the purposes of this article, but we did
not want the reader to think that projective geometry exhausts its role with the study of points and lines. We should
add that, like it often happens in mathematics, the theory of projective geometry and its developments has taken a life
of its own, and it is now one of the most fertile and successful fields in all of mathematics.}

{\selectlanguage{english}
To begin our analysis of the evolution of the prospective point of view in painting, we will look at a few paintings
from the early renaissance. Specifically, in section 2, we will look at two Tuscan painters: Giotto, whose worldwide
fame rests on his fresco cycle in Padova (where he depicted the life of Jesus and the life of the Virgin Mary), and
possibly (attribution is disputed) on his frescos in Assisi (where he depicted the life of San Francesco), and the
equally important Duccio di Buoninsegna, whose \textit{Maestà} is visible at the Duomo in Siena.}

{Giotto was
considered, at the time, the greatest living painter, and he is usually credited with being the link between the
Byzantine style and the Renaissance, and the first to adopt a more naturalistic style. \textcolor[rgb]{0.2,0.4,1.0}{
}Giotto was an attentive observer of reality, as we can see by looking at the faces and figures in his paintings, but
because of the lack of appropriate technique, his approach to architecture appears a mixture of artificial and
fantastic.\footnotemark{} In this section, we consider some of his works, as well as Duccio's paintings, to highlight
both their early understanding of the need for new ideas, as well as their insufficient clarity on what those ideas
would need to be.}

\footnotetext{\textrm{The reader is referred to [18], [26] for a careful reconstruction of the path from natural to
artificial perspective in the Middle-Ages and Renaissance. See also the bibliographies [19], [21], [25].}}

{Section 3 is devoted to the mathematical description of the
process that is necessary for a faithful representation of a three-dimensional scene on a canvas. This section is where
we are able to introduce the basic ideas that will lead to the projective space. How Leon Battista Alberti and Piero
della Francesca understood such ideas is the subject of Section 4, where we go back to the original texts, and
paintings, to illustrate the way in which the theory of projective geometry was applied in these more advanced works
from the Renaissance. To be precise, we will show that in fact Leon Battista Alberti did not fully justify his
technique (\textit{costruzione legittima}), and so we have an example of a process which seems to work, while its own
developers are not yet fully aware of its theoretical justification. The last Section, before our final conclusions,
inverts the process, so to speak. Instead of discussing how to use geometry to represent a scene on the canvas, we will
take a painting as a starting point, to reconstruct what the scene that the painter had in mind must have been. This is
an interesting exercise, not only for the mathematician, but for the art historian as well, since this reconstruction
can help shed light on some interpretation issues, as we will discuss in more detail in the section. }

\section{2. Early steps: Giotto (1267{}-1337) and Duccio di Buoninsegna  (1255/60{}-1318/19)}

{\selectlanguage{english}
If one takes a look at any of Giotto's frescos, the first thing that jumps to the eye is a really distorted sense of
distances, positions, and sizes of the elements of the pictorial composition. As we see below (figure 1) in a fresco
that represents San Francesco who chases away the devils from the city of Arezzo, the buildings have odd angles, the
figures are too big, and it looks like we are watching the scene both from the top and from the side (note how we see
the side of the walls surrounding Arezzo, but also the buildings inside the walls themselves). What is going on?}

{\selectlanguage{english}
The answer to this question lies in the fact that Giotto is one of those painters who found themselves in a moment of
epochal transformation. A moment in which painters understood that the way we see objects, and the way objects are, do
not coincide. More specifically, when we think of a table, when we touch a table, we deal with a rectangle. This is
what most tables are, and if we close our eyes and simply touch the table, we perceive a rectangle. Opposite sides are
parallel, and the angles between contiguous sides are right angles (ninety degrees). But when we look at a table, or
when we try to draw a table, something completely different appears. Now the angles become acute or obtuse, depending
on where we are looking, and the parallel sides may not appear parallel anymore.}

\bigskip

{\centering
\includegraphics[width=2.8799in,height=3.3335in]{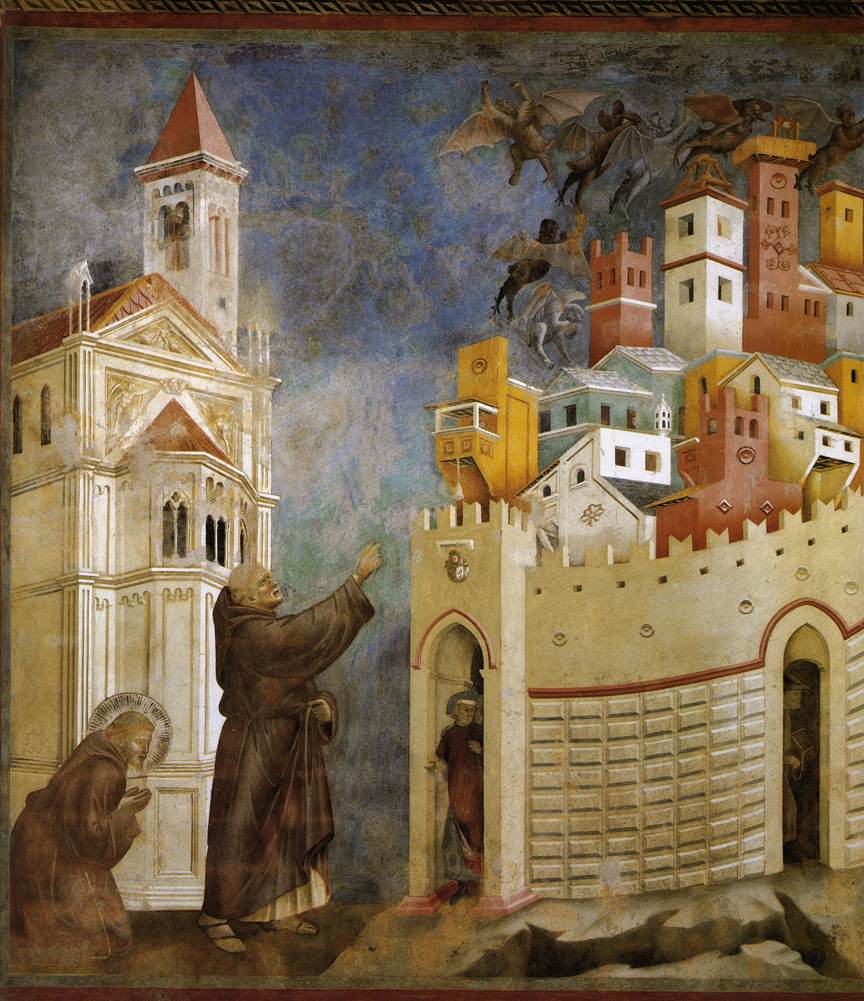}
\par}
{\centering\selectlanguage{english}
\foreignlanguage{italian}{Figure 1. Giotto, }\foreignlanguage{italian}{\textit{La cacciata dei diavoli da
Arezzo}}\foreignlanguage{italian}{, scene from ``Storie di San Francesco'', (1295-1299), fresco, Basilica Superiore di
Assisi.}
\par}

\bigskip

{\selectlanguage{english}
So, the painter has to recognize a complex shift: if the table has to look right, it has to be drawn wrong. Instead of a
rectangle, something else has to be drawn, in order to trick the viewer's brain into recognizing a properly positioned
table. If the painter were a mathematician, he would recognize that there are two geometries that conflict with each
other: the geometry of touching (the geometry of sculpture), and the geometry of seeing (the geometry of painting). But
this must have been incredibly difficult for Giotto and his contemporaries back in the fourteenth century. This
difficulty explains why his frescos appear so odd, and why the angles in the buildings that he depicts are so
un-lifelike. It is because Giotto understood that right angles do not always appear as right, but in fact they need to
be depicted as acute or obtuse. But he did not grasp, for example, the fact that parallel lines don't always appear
parallel. The unrealistic sizes of the figures in his frescos are a consequence of a similar cognitive dissonance.
Giotto realized that objects that are closer to us appear larger than objects at a distance. But he lacked the
mathematics to figure out the precise proportions that should be used. As we will see in Section 4, it will only be
with Leon Battista Alberti (1404-1472) that a precise method to address this issue will be developed.}

{\selectlanguage{english}
This conflict is quite apparent in another great contemporary of Giotto, namely Duccio di Buoninsegna. In his
\textit{Maest}à, there is a section where Duccio paints a \textit{Last Supper }(figure 2). }

\bigskip

{\centering 
\includegraphics[width=2.9736in,height=2.6957in]{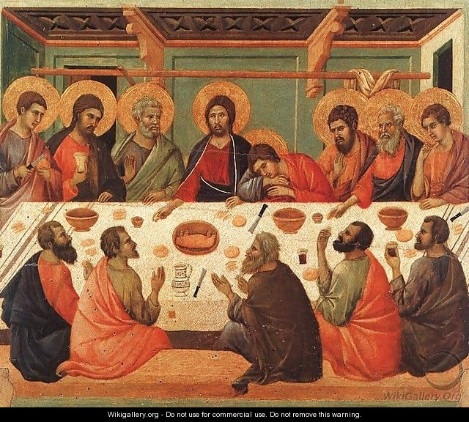}
\par}
{\centering\selectlanguage{english}
\foreignlanguage{italian}{Figure 2. Duccio di Buoninsegna, }\foreignlanguage{italian}{\textit{L'ultima
Cena}}\foreignlanguage{italian}{, scene from the back of the ``Maestà'', (1308-11), tempera on wood, Museo dell'Opera
del Duomo, Siena}
\par}
\bigskip
{\selectlanguage{english}
The central object in any such painting is the table, and when we look at this representation, we have the impression
that the plates on the table are on the verge of falling on the floor. The reason for such an impression is that the
table is represented not as a rectangle (Duccio like Giotto realized that this would not have worked), nor as a
trapezoid (which is the correct representation). Rather, it is a parallelogram, in which the right angles are
eliminated (as they should), but the parallelism among sides is preserved, thus offering a totally inadequate
representation. One should also look at the ceiling and the beams in the ceiling itself. In the room, such beams are
clearly parallel, and we know (we will get into more details later) that parallel lines must be represented as
converging to a point. But, as we see in the modified picture below, while Duccio understands this, he seems not to
know that all lines parallel to each other must converge in the same point. Instead, as we see (figure 3), the internal
beams converge on the figure of Christ, while the external beams converge on the table. The outcome is a ceiling that
is clearly wrong.}

\bigskip

{\centering 
\includegraphics[width=3.1563in,height=2.8173in]{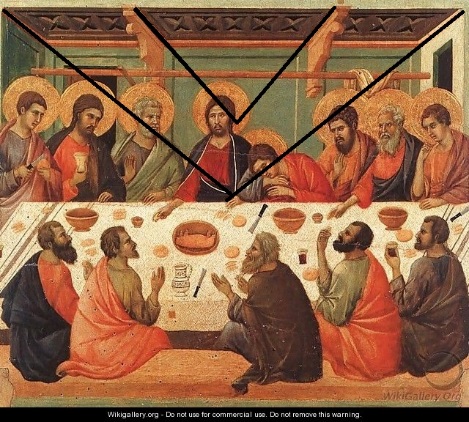}
\par}
{\centering\selectlanguage{english}
Figure 3
\par}
\bigskip
\bigskip

{\selectlanguage{english}
These two examples are not offered to demean these great painters, but rather to suggest how complex it must have been
for those living in the XIII and XIV century, to realize how to go from Euclidean geometry to projective geometry. As
we will see in the next sections, this process will lead to the understanding that a fundamental mathematical truth is
hidden somewhere. And it was because of a few artists with great mathematical background, that this was finally
understood. \ Before we get to that point, however, we will take a brief mathematical detour.}

\section{3. The mathematics of perspective}

{\selectlanguage{english}
It is an interesting challenge to illustrate in modern terms how the effort to understand vision and painting leads to a
new geometry, that mathematicians call projective geometry. \ }

{\selectlanguage{english}
First of all, notice that the main approach of a painter to this problem does not discuss how the eye and the brain
allow us to see the surrounding world, but only how the light and the colors reach the eye: in fact this is the
environment in which a painter can mainly intervene. It is therefore reasonable to think the eye as a point of our
3-dimensional space, and set its position as the origin $O$ of a system of Cartesian coordinates $(x,y,z)$.
\footnote{\textrm{ This is clearly a simplification that does not model the anatomical aspects of vision.}}}

{\selectlanguage{english}
We can then base our study on the experimental fact that light rays essentially propagate along straight lines in space,
and hence assume that the vision relies upon what all the infinite rays entering $O$ bring to the
eye.\footnote{\textrm{ This interpretation of vision, a revolutionary one in the Renaissance, was due mainly to ibn
Al-Haytham. also known as Alhazen, a well-known authority of the 11th century. Indeed, this visual theory was based on
his }\textrm{\textit{Book of Optics}}\textrm{ (}\textrm{\textit{Kitab al-Manazir}}\textrm{) [3]. On this topic see [5],
[6], [7].}} Every ray entering $O$ brings a colored point, which comes from the object being viewed (possibly
the sky). \ Therefore each ray entering the origin contributes to the vision with a colored point. Of course $O$
brings no contribution to the vision. }

{\selectlanguage{english}
Let us now imagine to be able to insert, between the observer, and the object which is observed, a canvas, possibly a
transparent one. Then, it is clear that the ray that joins the object to the observer will intersect the canvas in one
point and one alone. That point, with the color that the object has, becomes like a pixel on the canvas, and the
entirety of the pixels that are generated in this way, is a faithful representation of the object itself. Note that if
we were to take two different canvases, the eye would not be able to distinguish a difference in the images (we hope
the reader will forgive our use of relatively modern terms such as pixel).}

{\selectlanguage{english}
The beautifully simple, but not at all easy, mathematical idea that we have just described can be expressed by saying
that we have represented each spatial ray $r$ entering the origin $O$ by means of one of its points only,
$P(r)$, that lies on the chosen canvas and contains all the information that the ray brings to the eye. Of
course, if we move the representing point $P(r)$ (with all the information that it carries) along the ray
$r$, (in other words we change the canvas) the resulting view will not be affected at all. \ In principle, the
representing point of a ray can be chosen to be \textit{any }point of the ray. In practice, it is usually chosen to
stay on a plane - the plane of the painting, its canvas -- or, in different, more mathematical contexts, on a sphere
(we will not discuss this more complex type of representation).}

{\selectlanguage{english}
To help the reader understand what follows, we suggest a simple experiment. As you go through the next few lines, we
would ask that you sit in front of a window, looking at whatever lies in front of you.\textstyleFootnoteSymbol{
}\footnote{\textrm{ The metaphor of the window was used by Leon Battista Alberti to explain his
}\textrm{\textit{costruzione legittima. }}\textrm{According to Gerard Simon, without the new ideas of ibn Al-Haytham on
vision, Alberti's window would not have been thinkable: one of the many examples of historical encounters between
Western and Arab cultures [20]. }}\ \ And now, as you sit, imagine your eye to be the origin $O$ of a system of
coordinates (figure 4). The x-axis of the system will exit from your eye (the origin) and points to your right, the
$y$-axis exits from $O$ and points forward towards the window, and finally the $z$-axis is the vertical line from
$O$ up. Using the Cartesian coordinates so established, we will call ${\pi}$ the plane of the window (we assume
you are sitting upright, and therefore the window is perpendicular to the $y$-axis). If we assume the distance from the
reader to the window to be one unit, we would mathematically express the equation of this plan as $y=1$ (the mathematical
notations will be useful in the sequel when we will write the equations of the transformations, but are not necessary
for the understanding of the basic ideas). We will also identify the ceiling of the room with the plane of equation
$z=1$, i.e. with the horizontal plane located at a distance of $1$ unit from the eye, above the head of the reader.}

\includegraphics[width=5.24in,height=3.65in]{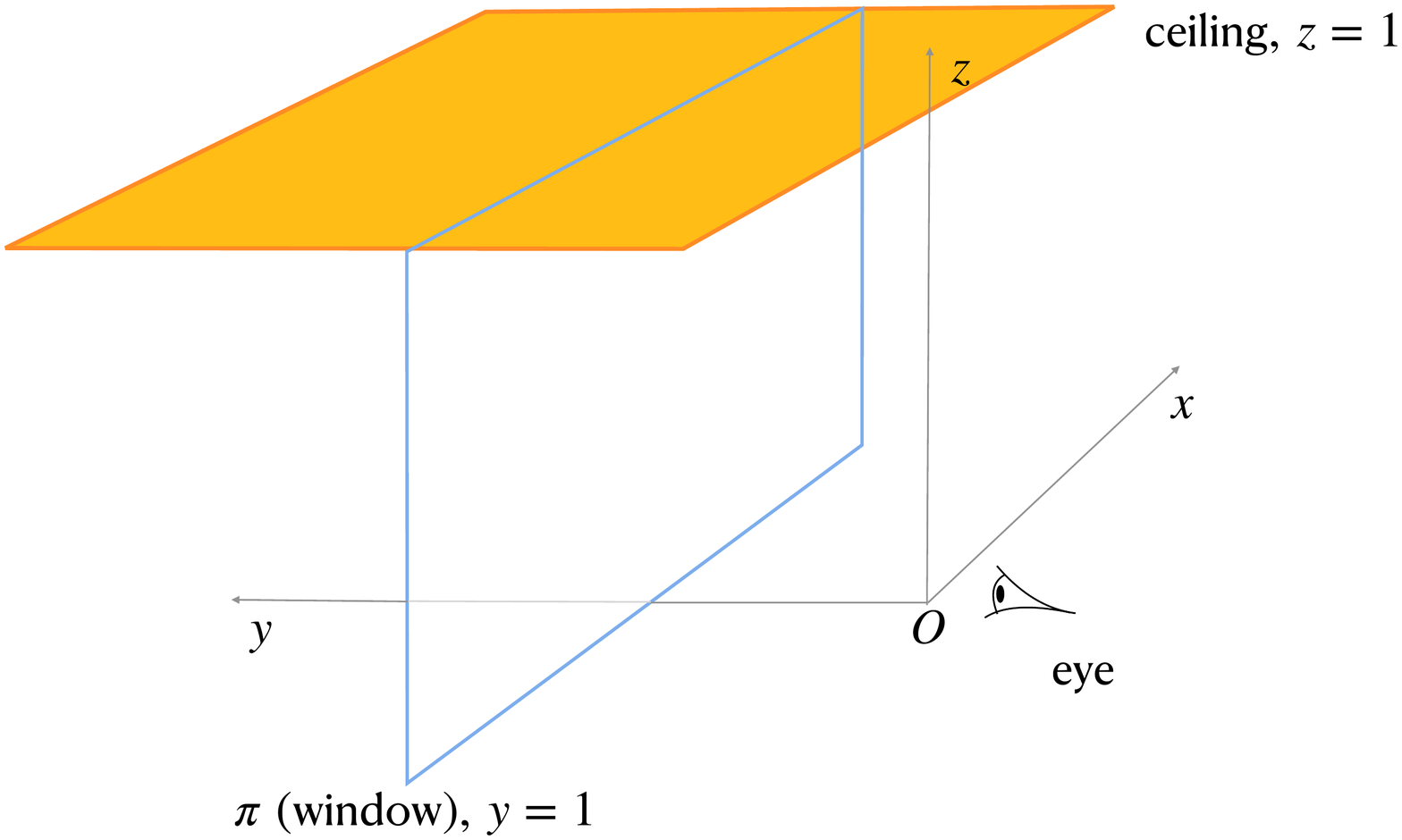}
{\centering\selectlanguage{english}
Figure 4
\par}
\bigskip
\bigskip
{\selectlanguage{english}
Now, two points of Cartesian coordinates $(x,y,z)$ and $(u,v,w)$ give the same contribution (pixel) $P(r)$ to the
vision if they belong to the same ray $r$ entering the origin (where the eye is located)\footnote{\textrm{
Mathematically, this means that there is a nonzero real number $t$ such that $(u,v,w)=t(x,y,z)=(tx, ty, tz)$. \ }}.
Therefore we will denote the contribution (pixel) $P(r)$ to the vision, given by the ray $r$ containing
the point $(x,y,z)$ and entering the origin, with the symbol $[x,y,z]$, \ and establish that $[x,y,z]=[tx,ty,tz]$ for all
nonzero real numbers t. \ The idea is that $[x,y,z]$ and $[tx,ty,tz]$ will indicate the same pixel, positioned on different
canvas.}

\includegraphics[width=5.24in,height=3.65in]{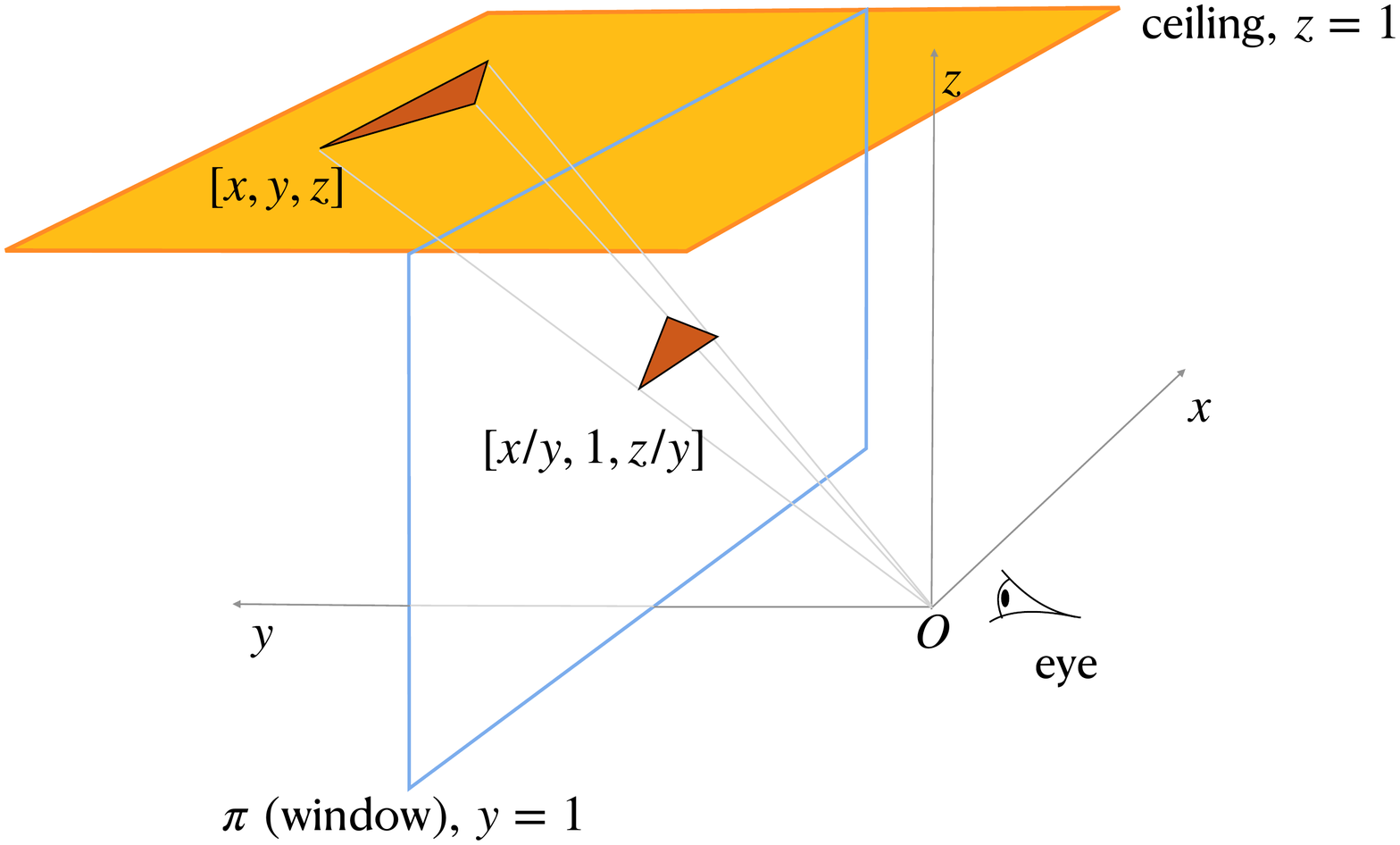}

{\centering\selectlanguage{english}
Figure 5
\par}
\bigskip
\bigskip

{\selectlanguage{english}
Since we have taken the window to be represented by the equation $y=1$, the pixel $[x,y,z]$ generated by $(x,y,z)$ will
correspond, on the window, to a point with $y=1$; this can only be obtained by taking $t=1/y$ and therefore the coordinate
of the pixel on the window will be $(x/y,1,z/y)$. (figure 5).}

{\selectlanguage{english}
Rays that enter the eye at $O$ will intersect the plane ${\pi}$ (just like when you are watching the countryside
from inside your home, the entering rays would all intersect the glass of the big window). If for all rays $r$
we place the representative point $P(r)$ on the plane ${\pi}$ \ than we have made the perfect theoretical
painting that represents the landscape the eye is watching.}

{\selectlanguage{english}
And here are a few surprises. A straight line $s$ of the observed landscape will be seen by the eye through all the rays
of the plane L that contains $O$ and the line $s$. We can then represent s on the painting ${\pi}$ as the
intersection of L and ${\pi}$. \ This demonstrates (in an empirical way) one of the first results of projective
geometry, namely the fact that a line is transformed (by projections) into another line (the reader is invited to
reflect on what would happen, however, if the line were one of the rays).}

{\selectlanguage{english}
With this in mind, if we are given the equations of a few beams of the ceiling of our room in the $3$-space, we can for
example compute the equations of the beams, and how they will appear in the painting ${\pi}$. As we have seen in
Duccio's example in the previous section, the issue of representing ceiling beams was in fact one of the most difficult
to understand.}

In our Cartesian environment, let us consider $3$ parallel beams in the ceiling $z=1$, of equations 
\begin{eqnarray*}
&&x=-1 \quad \hbox{and} \quad z=1\\
&&x=0 \quad \hbox{and} \quad z=1\\
&&x=1 \quad \hbox{and} \quad z= 1,
\end{eqnarray*}

{\selectlanguage{english}
respectively\footnote{\textrm{ The reader will note that we use two equations to represent a line. This is because
}\textrm{{one}}\textrm{ can think of a line as the intersection of two planes, each one with its own
equation. In this particular case, we are looking at lines which are on the ceiling (and so all of their points have
$z=1$) , but also that are perpendicular to the x-axis and therefore have a fixed value for $x$ (in the three cases,
respectively, $x=1, x=0, x=-1$).}}. These three beams are all parallel to the y-axis as indicated in figure 6.}

\bigskip

{\centering 
\includegraphics[width=5.24in,height=3.6445in]{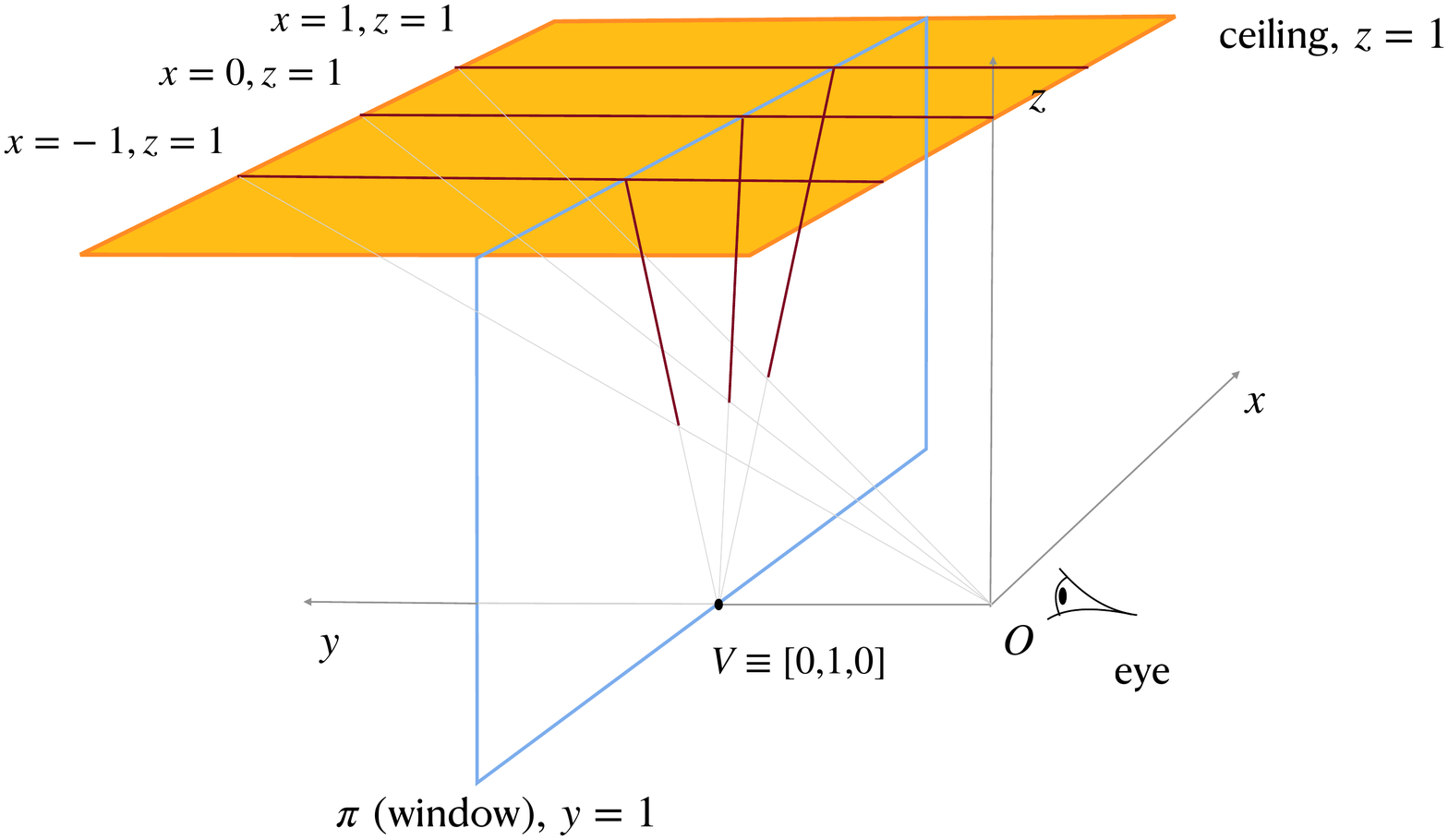}
\par}
{\centering\selectlanguage{english}
Figure 6
\par}
\bigskip
\bigskip

{\selectlanguage{english}
Note that a point on the first beam (the one with equations $x=-1, z=1$) will have coordinates $(-1, y, 1)$, where the $x$ and
the $z$ coordinates are fixed because of the planes and the $y$-coordinate is free to range in any way.}

{\selectlanguage{english}
If we now search for the three sets of points that represent the contributions to the vision (pixels), coming from the
three beams, we get (with arbitrary $y$)}
\[
[-1,y,1]
\]
\[
[0,y,1]
\]
\[
[1,y,1].
\]

As we have already pointed out, to place them on the painting ${\pi}$ of equation $y=1$, we just divide the (so called
homogeneous) coordinates by (the arbitrary nonzero) $y$ and get, with the established notations:

\begin{eqnarray}\label{(1)}
&&[-1/y, 1, 1/y]\\
&&[0, 1, 1/y]\nonumber  \\
&&[1/y, 1, 1/y],\nonumber
\end{eqnarray}

i.e., by putting $u=1/y$,
\bigskip
\begin{eqnarray}\label{(2)}
&&[u,1, -u]\\
&&[0,1, u]\nonumber \\
&&[u, 1, u].\nonumber
\end{eqnarray}

{\selectlanguage{english}
These are the equations of three straight lines in ${\pi}$ that are not only nonparallel, but that all meet at the point
$V$=(0,1,0).\footnote{ Strictly speaking, the point $V$ cannot be achieved because u is never zero, but
we think it is clear what we are describing here.} (figure 6).}

{\selectlanguage{english}
If the reader has followed the process, (s)he should be noticing that this process establishes a certain correspondence
between the ceiling and the canvas. Every point of the ceiling has a corresponding point on the canvas, but not every
point on the canvas comes from a point on the ceiling. \ Indeed, it is apparent that the points of the canvas y=1 with
a negative z{\textless}0 coordinate cannot come from the ceiling: the reader will immediately see that these points
come from the floor (ground) of the observed landscape (floor of equation, say, z=-1). Well, now all seems to be well
understood{\dots}, but where is the point $V$ coming from? If one tries to reconstruct the process we just
described, one will realize that in fact the point $V$ does not come either from any point on the ceiling, or
from any point of the floor and this throws a monkey-wrench in our construction. What can be done to fix this apparent
irregularity? What is the meaning of this surprising difficulty?}

{\selectlanguage{english}
If we analyze carefully what we have done so far, we will notice that the point $V$ is approached\textbf{ }on the
painting\textbf{ }${\pi}$ by the pixels contributed by the rays coming from points of any of the three beams very far
from the origin: when y becomes arbitrarily big in (1), $1/y= u$ approaches $0$ (see (2)). \ If we now think of a few lines
of the floor, parallel to the y axis, and repeat for the floor $z=-1$ the procedure used for the ceiling, we will see
that the point $V$is approached\textbf{ }on the painting\textbf{ }${\pi}$ by the pixels contributed by the rays
coming from points of any of these lines of the floor, very far from the origin.}

{\selectlanguage{english}
In some sense, the point $V$ (which in the painter terminology is called \textit{the vanishing point of the
painting} \footnote{\textrm{ Lucretius, in his }\textrm{\textit{De rerum natura}}\textrm{ describes the vision of a
colonnade which extends in front of us \ in the passage we used as incipit to this article: «Again, a colonnade may be
of equal line from end to end and supported by columns of equal height throughout, yet, when its whole length is
surveyed from one end, it gradually contracts into the point of a narrowing cone, completely joining roof to floor and
right to left, until it has gathered all into the vanishing point of the cone.» Lucretius.~}\textrm{\textit{On the
Nature of Things, }}\textrm{[17,}\textrm{\textit{ }}\textrm{IV, pp. 426-431].}\textrm{\textit{ }}\textrm{This
fascinating piece of poetry constitutes the first description of the vanishing point that has reached us. }}) is the
image of the point on the ceiling that would belong to each of the beams, if they could continue to
infinity.\footnote{\textrm{ For an extensive explanation, see e.g. [11].}} In the same way, the point $V$ is the
image of the point on the floor that would belong to each of the chosen parallel lines, if they could continue to
infinity. \ But of course the three beams, and the chosen lines of the floor, are parallel, and they have no point in
common. What is happening? The answer (whose mathematical formalization we will describe shortly) is that in order for
the correspondence between the canvas and the system ceiling-floor to be complete, we need to add a point (which
doesn't belong either to the ceiling or to the floor), which is the intersection both of the parallel beams of the
ceiling and of the chosen parallel lines of the floor. In fact, as we will discover shortly, even this addition will
not be enough. Indeed, we will need to add to the system an entire line, in order to reconstruct a perfect
correspondence.}

{\selectlanguage{english}
To understand this last point, consider now a family of parallel beams on the ceiling that are not parallel to the $y$
axis. Consider for instance the three beams of equations (figure 7)}
\begin{eqnarray}\label{(3)}
&&x=y \quad \hbox{and}\quad z=1\\
&&x=y-1 \quad \hbox{and}\quad z=1\nonumber \\
&&x=y-2 \quad \hbox{and}\quad  z=1.\nonumber
\end{eqnarray}

{\selectlanguage{english}
These three beams contribute to the vision with the pixels denoted by, for arbitrary $y$}

\begin{eqnarray}\label{(4)}
&&[y, y, 1]\\
&&[y-1, y, 1]\nonumber \\
&&[y-2, y, 1],\nonumber
\end{eqnarray}

{\selectlanguage{english}
which placed on the painting ${\pi}$ \ (of equation $y=1$) become, for arbitrary nonzero $y$}

\begin{eqnarray}\label{(5)}
&&[y/y, 1, 1/y] = [1, 1, 1/y]\\
&&[(y-1)/y, 1, 1/y] = [1-1/y,1, 1/y]\nonumber \\
&&[(y-2)/y, 1, 1/y] = [1-2/y, 1,1/y],\nonumber
\end{eqnarray}

{\selectlanguage{english}
i.e., for an arbitrary nonzero $u$}

\begin{eqnarray}\label{(6)}
&&[1, 1, u]\\
&&[1-u, 1, u]\nonumber \\
&&[1-2u, 1, u].\nonumber
\end{eqnarray}

{\selectlanguage{english}
These are the equations of three straight lines in ${\pi}$ that, again, are not only nonparallel, but that all meet at
the point $W$= $(1,1,0)$ (figure 7). Again, the point $W$ does not contribute with a pixel coming from the
ceiling, and anyway belongs to the painting ${\pi}$. \ $W$ is called \textit{vanishing point for the given
family of parallel beams.} \ The point $W$is approached\textbf{ }on the painting\textbf{ }${\pi}$ by the pixels
contributed by the rays coming from points of any of the three beams listed in (3), very far from the origin: when y
becomes arbitrarily big in (5), $1/y= u$ approaches 0 in (6). If we now think of a few lines of the floor, obtained by
substituting $z=-1$ in place of $z=1$ in formulas (3), and repeat the procedure used for the ceiling, we will see that
the point $W$is approached\textbf{ }on the painting\textbf{ }${\pi}$ by the pixels contributed by the rays
coming from points of any of these lines of the floor, very far from the origin.}

{\selectlanguage{english}
The family of parallel lines that we considered in (3) are actually parallel to the bisecting line of the first and
third quadrant of the $xy$ plane of equation $z=1$ (the ceiling): these lines could be diagonals of a square tessellation
of the ceiling. In this situation, we see that the distance between the vanishing point of the painting
$V$=$(0,1,0)$ and the vanishing point $W$=(1,1,0) coincides with the distance of the eye of the observer
from the plane of the painting ${\pi}$ (!!). \ The distance of the eye of a painter from his painting can be encoded
in the painting itself.  This is the reason why Leon Battista Alberti e Piero della Francesca called $W$ the
\textit{distance point}. Notice that the projective approach made the identification of the point $W$
immediate.}

\bigskip

\includegraphics[width=5.24in,height=3.65in]{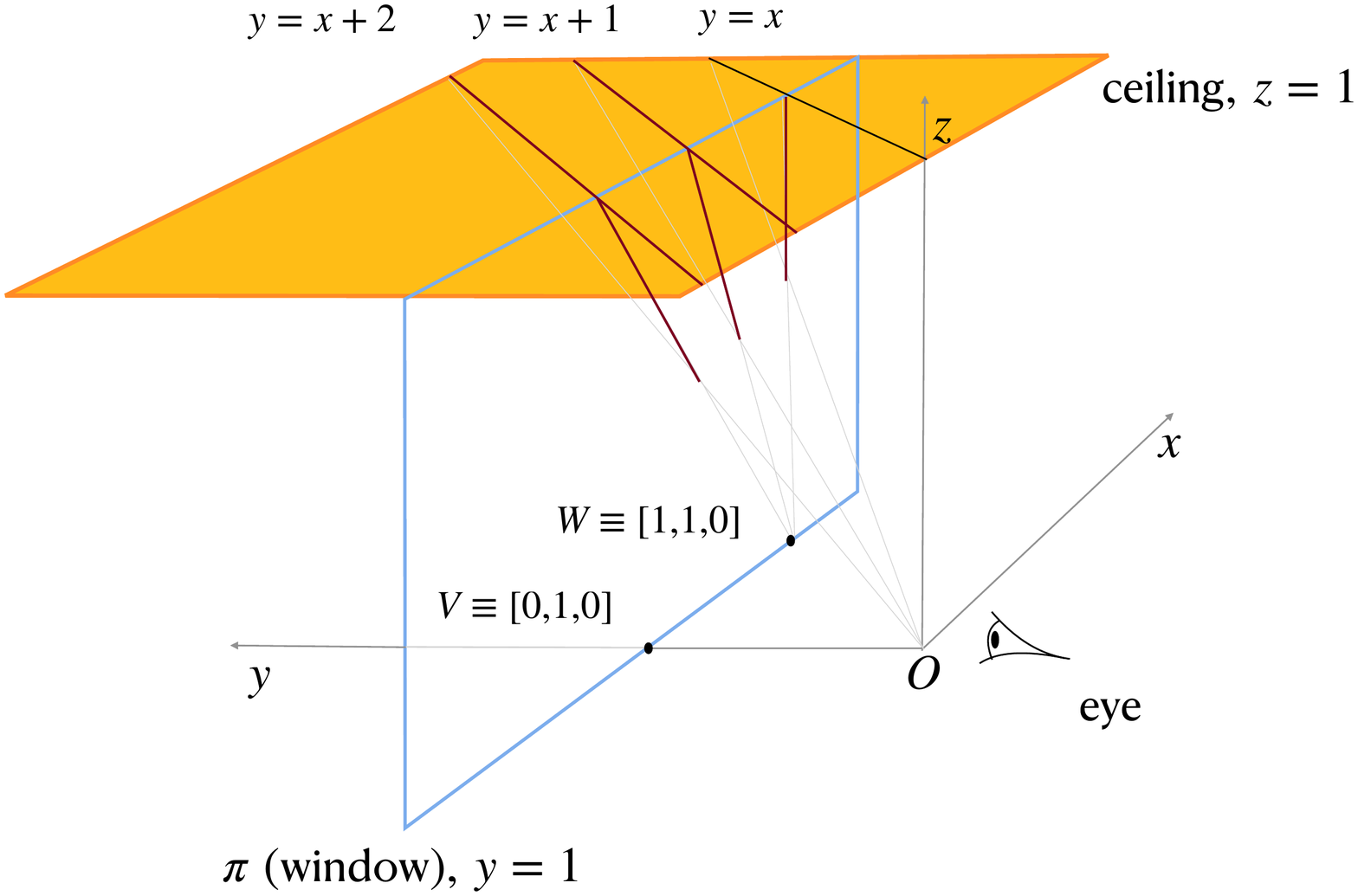}

{\centering\selectlanguage{english}
Figure 7
\par}

\bigskip

{\selectlanguage{english}
The phenomenon we described above is not limited to a particular family of parallel lines. More generally, we will show
that every family of parallel lines on the ceiling is represented, on the canvas, by a family of lines that converge to
a point that doesn't come either from a point of the ceiling, or from a point of the floor, and that needs to be added
in order to complete the correspondence of the canvas with the system ceiling-floor. And all these new points that we
will add (and which we call \textit{improper points} or \textit{points at infinity}), will eventually be on a line
(\textit{improper line} or \textit{line at infinity}), whose pictorial meaning we will describe in the next few pages.}

{\selectlanguage{english}
Let us therefore consider an arbitrary family of parallel lines in the ceiling of equation $z=1$. Three beams from this
family have equations, for any nonzero $m$, }

\begin{eqnarray}\label{(7)}
&&x=(y-1)/m \quad \hbox{and} \quad z=1\\
&&x=(y-2)/m \quad \hbox{and} \quad  z=1\nonumber \\
&&x=(y-3)/m \quad \hbox{and} \quad  z=1.\nonumber
\end{eqnarray}

{\selectlanguage{english}
These three beams contribute to the vision with the pixels denoted by, for arbitrary $y$}

\begin{eqnarray}\label{(8)}
&&[(y-1)/m, y, 1]\\
&&[(y-2)/m, y, 1]\nonumber \\
&&[(y-3)/m, y, 1],\nonumber
\end{eqnarray}
{\selectlanguage{english}
which on the painting ${\pi}$ \ (of equation $y=1$) become, for arbitrary nonzero y}

\begin{eqnarray}\label{(9)}
&&[(y-1)/my, 1, 1/y] = [1/m-1/my, 1, 1/y]\\
&& [(y-2)/my, 1, 1/y] = [1/m-2/my,1, 1/y]\nonumber \\
&&[(y-3)/my, 1, 1/y] = [1/m-3/my, 1,1/y],\nonumber
\end{eqnarray}

{\selectlanguage{english}
i.e., for an arbitrary nonzero $u$}$U$

\begin{eqnarray}\label{(10)}
&&[1/m-u, 1, u]\\
&&[1/m-2u, 1, u]\nonumber \\
&&[1/m-3u, 1, u].\nonumber
\end{eqnarray}

{\selectlanguage{english}
These are the equations of three straight lines in ${\pi}$ that, again, are not only nonparallel, but that all meet at
the point $U= (1/m,1,0)$. Again, the point $U$ does not co$U$ntribute with a pixel coming from the ceiling,
and anyway belongs to the painting ${\pi}$. The point $U$ is called \textit{vanishing point for the given family of
parallel lines}. \ If we now consider on the floor the family of lines analogous to the family described in (7), but
with $z=-1$, we still find the point $U$, which\textbf{ }does not contribute with a pixel coming from the floor,
and anyway belongs to the painting ${\pi}$.}

{\selectlanguage{english}
It is clear now that (being m arbitrary) the collection of all vanishing points of families of parallel lines of the
ceiling or of the floor, contribute to the vision with all the pixels that on the painting ${\pi}$ of equation $y=1$ are
of the form}
\[
[u,1,0],
\]

{\selectlanguage{english}
i.e. with the line of the painting ${\pi}$ of equation ($y=1$ and) $z=0$. For obvious and charming reasons, this line is
called the \textit{horizon}! (figure 8).}

\includegraphics[width=5.24in,height=3.65in]{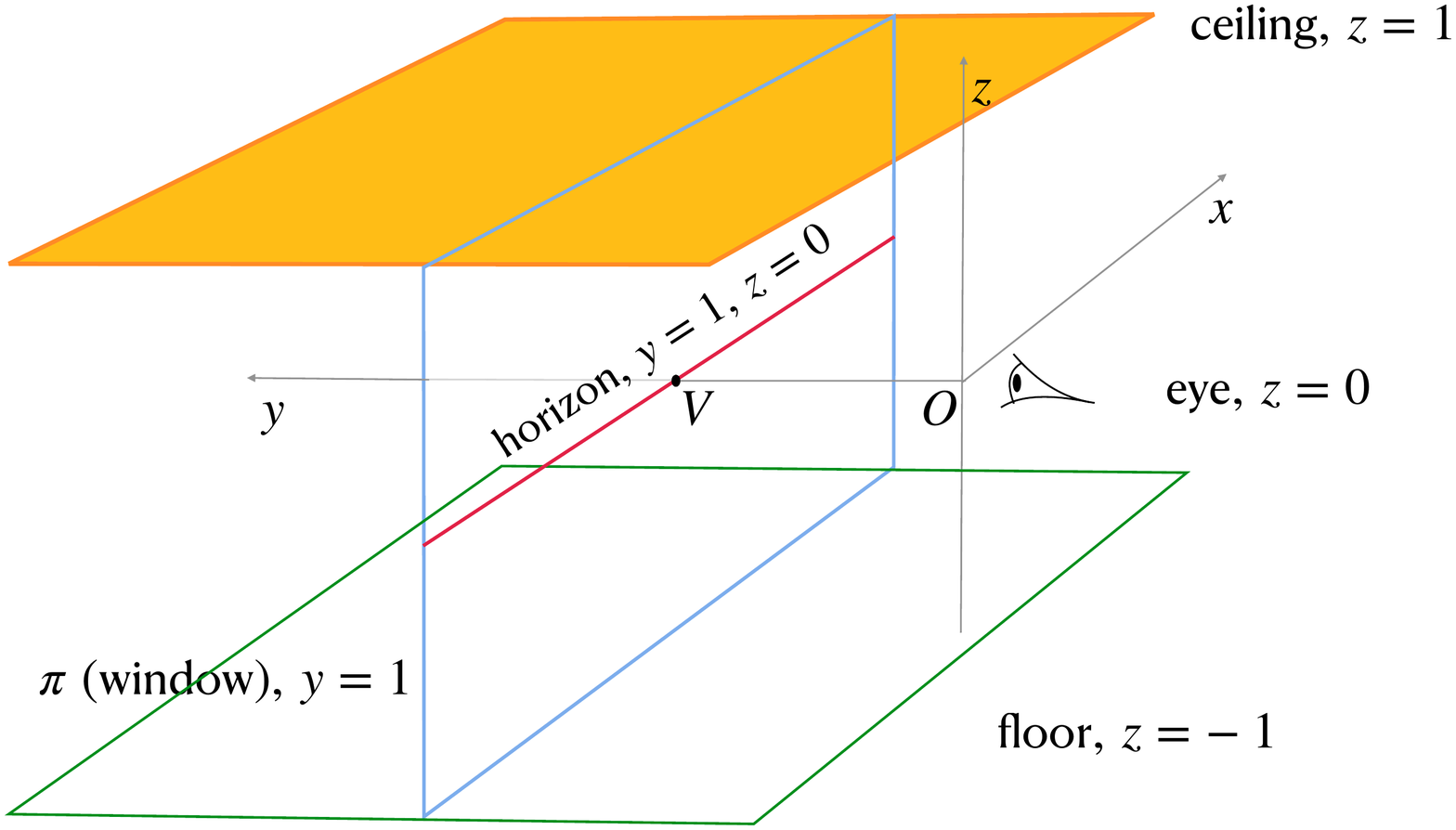}

{\centering\selectlanguage{english}
Figure 8
\par}
\bigskip

{\selectlanguage{english}
This construction is beautifully illustrated in the following painting of Andrea del Castagno (figure 9), in which we
have highlighted (figure 10) three families of parallel lines, with the corresponding vanishing points and the
resulting horizon (the reader is advised not to highlight such lines when visiting an art museum!). }

\bigskip

\includegraphics[width=5.461in,height=2.2264in]{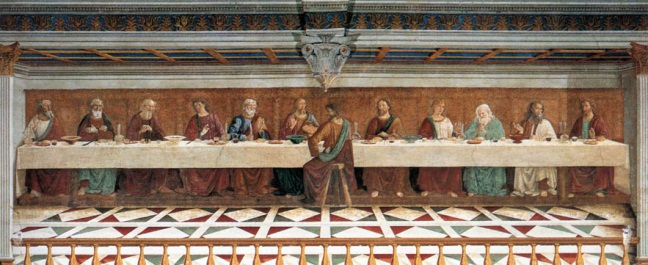}

{\centering\selectlanguage{english}
\foreignlanguage{italian}{Figure 9. Domenico Ghirlandaio, }\foreignlanguage{italian}{\textit{Ultima
cena}}\foreignlanguage{italian}{, (\~{} 1476), affresco, Cenacolo della Badia di Passignano, Abbazia di San Michele
Arcangelo a Passignano, Tavarnelle Val di Pesa, Firenze,}
\par}

\bigskip
\bigskip

\includegraphics[width=5.5827in,height=2.2264in]{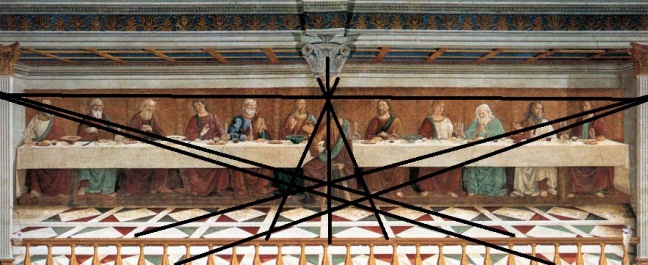}

{\centering\selectlanguage{english}
Figure 10
\par}

\bigskip

{\selectlanguage{english}
In order for the correspondence that we have described to hold for every point, we need to add an improper point for
every direction of lines. So, we now have an entire line of improper points on the system ceiling-floor, usually called
the improper line. The image of the improper line under the correspondence we have described is thus the horizon on the
painting. }

{\selectlanguage{english}
In this way we have presented a formalized mathematical method to move the representative point of each ray of light
(originating at the ceiling or at the floor) to the painting ${\pi}$. }

{\selectlanguage{english}
The reader will appreciate the beautiful symmetry that is emerging. Just like the images of two parallel lines (whether
on the ceiling or on the floor) intersect in a point, the vanishing point, that we imagine to be the image of the
improper point shared by the parallel lines, so the images of two parallel planes (the ceiling and the floor) intersect
in a line, the horizon, that is the image of the improper line that ceiling and floor share. A marvelous symmetry
indeed!}

{\selectlanguage{english}
The attentive reader will, however, note that while we have added a line (an improper one) to the system ceiling-floor,
the full correspondence will require the addition of a line to the (infinite) canvas as well. Indeed, if we are now
trying to describe (on the canvas) the line which is represented on the ceiling z=1 by y=0, we easily see that this is
not possible. Pictorially, this is a consequence of the fact that the painter cannot represent, in the painting, the
points that are vertically above his head. Mathematically, this is a consequence of the fact that the plane y=0 does
not intersect the plane ${\pi}$ \ given by=1. Finally, if one looks at the coordinates (x,y,1) of a point on the
ceiling, and allows y to become zero, one obtains the point (x,0,1) which does not belong to the painting. Just like we
did before, we would now need to add all these points (for all values of x) to the canvas, and thus complete the plane
of the painting with an improper line.}

{\selectlanguage{english}
In this new correspondence the following happens:}

\liststyleWWNumv
\begin{itemize}
\item {\selectlanguage{english}
Every point in the ceiling and on the floor (except those with y=0) is represented by a point on the canvas.}
\item {\selectlanguage{english}
Every point on the canvas (except those on the horizon) are the representation of a point on the ceiling or on the
floor.}
\item {\selectlanguage{english}
The points of the horizon can be thought of as images of the improper line that we have added to the system
ceiling-floor.}
\item {\selectlanguage{english}
The points on the ceiling with y=0 are represented on the improper line that we have added to the canvas.}
\end{itemize}
\section{4. Leon Battista Alberti (1404{}-1472) and Piero della Francesca (1416/17{}-1492)}
{\selectlanguage{english}
Section 2 described some of the uncertainties that were plaguing the painters of the early Renaissance. Despite these
uncertainties, these painters were feeling the strong need to change the nature and the subjects of their work. The
interest was slowly shifting away from the ascetic body of the teachers of the medieval scholastics, and was turning to
three-dimensional figures, the divine \textit{Maestà }inside gothic churches, or the suggestive backgrounds of battles
where the powerful soldiers and the vigor of the horses could find an effective representation. A philosophical
development was forcing the painters towards a new understanding of their art as evidenced in the work of artists such
as Paolo Uccello (figure 13), Mantegna (figure 11), Masaccio, and the Giambellino (Giovanni Bellini) (figure 12).}
\bigskip

{\centering 
\includegraphics[width=3.4866in,height=2.9555in]{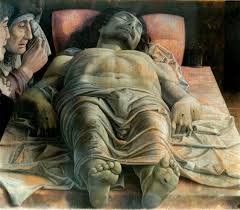}
\par}
{\centering\selectlanguage{english}
\foreignlanguage{italian}{Figure 11. Andrea Mantegna, }\foreignlanguage{italian}{\textit{Cristo morto
}}\foreignlanguage{italian}{(1475 {}- 1478), Tempera on canvas, Pinacoteca di Brera, Milano}
\par}

{\centering 
\includegraphics[width=3in,height=3.9307in]{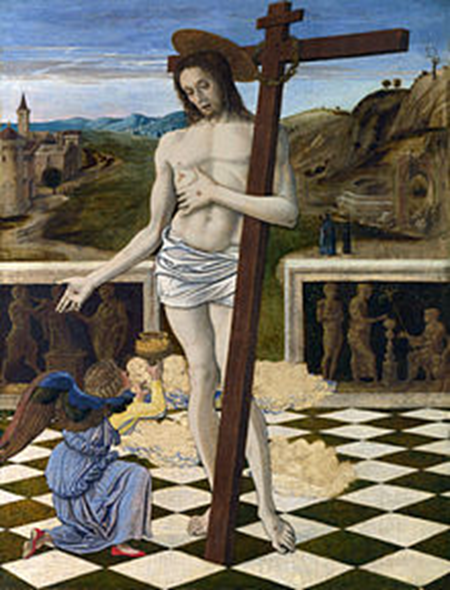}
\par}
{\centering\selectlanguage{english}
Figure 12. Giovanni Bellini, \textit{The Blood of the Redeemer}, (1460{}-1465), The National Gallery, London.
\par}

\bigskip

\bigskip

\section[\ ]{
\includegraphics[width=5.3134in,height=2.939in]{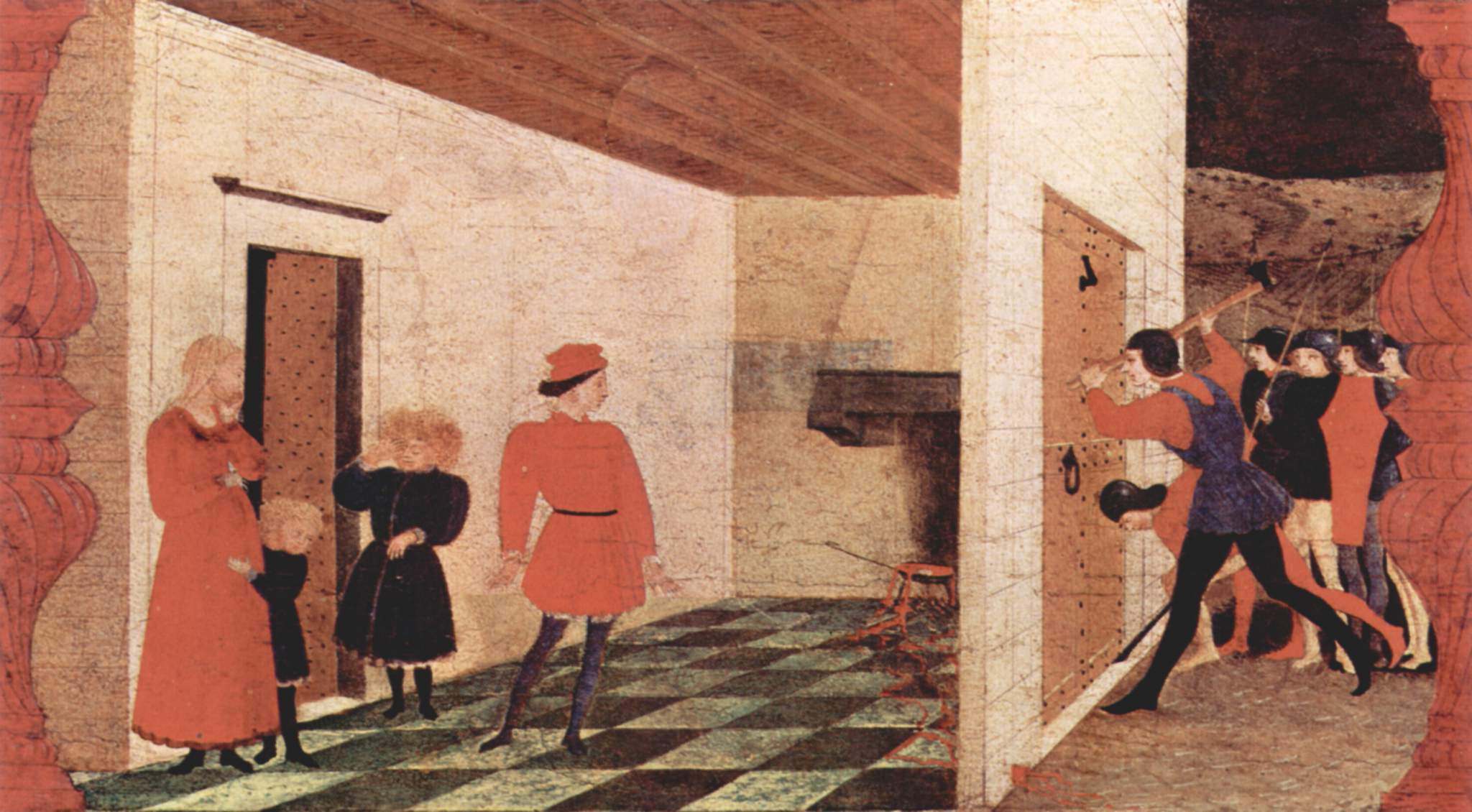}
\foreignlanguage{italian}{\ }}
{\centering {Figure 13. Paolo Uccello,
}\foreignlanguage{italian}{\textit{Predella del Miracolo dell'ostia profanata}}\foreignlanguage{italian}{, (1467-1468),
tempera on wood \ Galleria Nazionale delle Marche, Urbino.}\par}

\bigskip

{\selectlanguage{english}
Among them Leon Battista Alberti (who wrote in 1435 the treatise \textit{De pictura praestantissima} [1], where he offers a practical
guide to perspective drawing) and the great painter and mathematician Piero della Francesca, who built on his knowledge
of Euclid and Alberti, to write (towards the end of the XV century) \textit{De prospectiva pingendi} [15], probably the
ultimate text on prospective in painting.\footnote{\textrm{ Field's essay [10] - an extensive comparison of the
treatment of perspective in Alberti's }\textrm{\textit{De pictura praestantissima}}\textrm{ [1] and Piero della Francesca's
}\textrm{\textit{De prospectiva pingendi }}\textrm{[15]}\textrm{\textit{{}--}}\textrm{ contains a historically
contextualized presentation of the main mathematical tools on which the theory and practice of perspective (and the
very initial basis of projective geometry) relied upon. See also [8]}}}

\bigskip
\hfill \parbox[r]{12cm}{For theory, when a separated from practice, is generally of very little use; but when the two chance to come
together, there is nothing that is more helpful to our life, both because art becomes much richer and more perfect by
the aid of science, and because the counsels and the writings of learned craftsmen have in themselves greater efficacy
and greater credit than the words or works of those who know nothing but mere practice, whether they do it well or ill.
And that all this is true is seen manifestly in Leon Batista Alberti, who, having studied the Latin tongue, and having
given attention to architecture, to perspective, and to painting, left behind him books written in such a manner, that,
since not one of our modern craftsmen has been able to expound these matters in writing, although very many of them in
his own country have excelled him in working, it is generally believed; such is the influence of his writings over the
pens and speech of the learned; that he was superior to all those who were actually superior to him in work. ]{For
theory, when a separated from practice, is generally of very little use; but when the two chance to come together,
there is nothing that is more helpful to our life, both because art becomes much richer and more perfect by the aid of
science, and because the counsels and the writings of learned craftsmen have in themselves greater efficacy and greater
credit than the words or works of those who know nothing but mere practice, whether they do it well or ill. And that
all this is true is seen manifestly in Leon Batista Alberti, who, having studied the Latin tongue, and having given
attention to architecture, to perspective, and to painting, left behind him books written in such a manner, that, since
not one of our modern craftsmen has been able to expound these matters in writing, although very many of them in his
own country have excelled him in working, it is generally believed; such is the influence of his writings over the pens
and speech of the learned; that he was superior to all those who were actually superior to him in
work.}\textstyleRimandonotaapiiii{ }\footnotemark{}}
\footnotetext{\textrm{ }\textrm{\textit{Live of Leon Battista Alberti}}\textrm{ in: }\textrm{Vasari's
}\textrm{\textit{Lives of the Artists}}\textrm{, [23]. }\foreignlanguage{italian}{\textrm{See Vasari's
}}\foreignlanguage{italian}{\textrm{\textit{Le
vite}}}\foreignlanguage{italian}{\textrm{:}}\foreignlanguage{italian}{\textrm{
``}}\foreignlanguage{italian}{\textrm{{\dots} Non è cosa che più si convenga alla vita nostra, sì perché l'arte col
mezzo della scienza diventa molto più perfetta e più ricca, sì perché gli scritti et i consigli de' dotti artefici
hanno in sé molto maggiore efficacia et acquistansi maggior credito che le parole o l'opere di coloro che non sanno
altro che il semplice esercizio, o bene o male che essi lo facciano: ché invero leggendo le istorie e le favole et
intendendole, un capriccioso maestro megliora \ continovamente e fa le sue cose con più bontà e con maggiore
intelligenza che non fanno gli illetterati. \ E che questo sia il vero si vede manifestamente in
Leon}}\foreignlanguage{italian}{\textrm{\textbf{\textit{ }}}}\foreignlanguage{italian}{\textrm{Batista Alberti, il
quale, per avere atteso alla lingua latina e dato opera alla architettura, alla prospettiva et alla pittura, lasciò i
suoi libri scritti di maniera che, per non essere stato fra gli artefici moderni chi le abbia saputo distendere con la
scrittura, ancora che infiniti ne abbiamo avuti più eccellenti di lui nella pratica'',
}}\foreignlanguage{italian}{\textrm{[22, }}\foreignlanguage{italian}{\textrm{3
voll}}\foreignlanguage{italian}{\textrm{., vol. III, pp.
283-284}}\foreignlanguage{italian}{\textrm{].}}\foreignlanguage{italian}{\textrm{ See also [24, ``Leon Battista
Alberti'', pp.178-184]. \ }}}

\bigskip
{This is how Vasari, in his \textit{Vite} [22], (essentially a collection of biographies
of painters, sculptors, and architects) described the polyhedric character of Leon Battista Alberti, (1404-1472),
architect, mathematician, humanist, musician, who was born in Genova but split his life between the papal court of
Rome, and the courts of the Este in Ferrara, of the Malatesta in Rimini, of the Gonzaga in Modena, and belonged the
circle of the Florentine humanists. It is important to note that his reflections on painting and sculpture were not
simply the byproduct of his interest on painting techniques and on how to prospectively represent the human body, but
rather were consequence of a deeper intellectual research.}
{\selectlanguage{english}

In the opening of his treatise \textit{De pictura praestantissima }[1]\textit{, }Leon Battista Alberti explicitly declares that he is
not writing as mathematician, but as a painter. And in fact it is clear that the aim of his treatise is to provide
painters a practical guide to the use of perspective, instead of investigating and discussing in detail the theoretical
aspects and features of that subject, which -- as he explicitly states - was certainly quite difficult and not yet well
discussed by any author.\footnote{\foreignlanguage{italian}{\textrm{ Leon Battista Alberti,
}}\foreignlanguage{italian}{\textrm{\textit{The Architecture }}}\foreignlanguage{italian}{\textrm{[2, p. 241]. See also
Leon Battista Alberti, }}\foreignlanguage{italian}{\textrm{\textit{On Painting}}}\foreignlanguage{italian}{\textrm{ [2,
p. 37]. }}}}

\bigskip
{\selectlanguage{english}
\hfill \parbox[r]{12cm}{ But throughout these whole Treatise I must beg my Reader to take Notice, that I speak of these Things, not as a
Mathematician, but as a Painter; for the Mathematician considers the Nature and Forms of Things with the Mind only,
absolutely distinct from all Kind of Matter: whereas it being my Intention to set Things in a Manner before the Eyes,
it will be necessary for me to consider them in a Way less refined. And indeed I shall think I have done enough, if
Painters, when they read me, can gain some Information in this difficult Subject, which has not, as I know of, been
discussed hitherto by any Author.} }

\bigskip
{\selectlanguage{english}
In modern terms, we could say that what Leon Battista Alberti was doing was actually to scientifically present an
algorithm that any painter could use to correctly set up the basics of his paintings, from the point of view of the
perspective. \ More specifically, Alberti wanted to give a practical tool to correctly set up the floor, the vanishing
point, and the horizon of a painting (see Section 3). Once done this, it was easier for a painter to fill in the
painting with all the rest in a reasonably coherent form. This point of view explains also the reason why Alberti is in
fact teaching how to represent in perspective a ground floor with square tiles: even in the case of an eventually
uniform ground floor, the hidden presence of a fine square grid would help a lot the skilled painter to place objects
and human figures properly and proportionally in the table.}

{\selectlanguage{english}
It is of great interest to examine the three steps of the algorithm proposed by Leon Battista Alberti to paint a
square-tile floor in perspective. In his treatise \textit{De pictura praestantissima,} Alberti considers a square painting ${\pi}$
whose side measures six \textit{braccia fiorentine} (in modern terms approximately 348/354 cm). \ Since the established
standard height of a human figure for a painter in those years was three braccia fiorentine, in practice Alberti chose
to place:}

\liststyleWWNumi
\begin{itemize}
\item {\selectlanguage{english}
the horizontal basis of the painting on the floor;}
\item {\selectlanguage{english}
the point of view of the painter on the straight line orthogonal to the center of the painting ${\pi}$;}
\item {\selectlanguage{english}
the vanishing point at the center of the painting itself.}
\end{itemize}
{\selectlanguage{english}
One other datum is that the square tiles of the floor to be painted have two sides parallel, and two orthogonal, to the
painting $\pi $. Finally, here is Alberti's \ costruzione legittima.\footnote{\textrm{ The Renaissance texts on
perspective are normally didactic manuals, whose authors take it for granted that perspective is a
}\textrm{\textit{vera scientia}}\textrm{. Alberti mentions but does not give a proof for his
}\textrm{\textit{costruzione legittima}}\textrm{. In his article, Elkins [9] presents an annotated (incomplete) proof
of Alberti's }\textrm{\textit{costruzione,}}\textrm{ taken from two propositions of Piero's
}\textrm{\textit{Prospective pingendi}}\textrm{. }}}

{\selectlanguage{english}
\textit{Step 1. \ Design the projections of the ``orthogonal'' straight lines of the floor on the painting
}${\pi}$\textit{. }\ This step can be done formally as explained in Section 3. It can be practically performed as
follows: it is enough to join each intersection of a straight line of the floor with the basis of the painting with the
vanishing point (figure 14).}

\includegraphics[width=5.5134in,height=2.7654in]{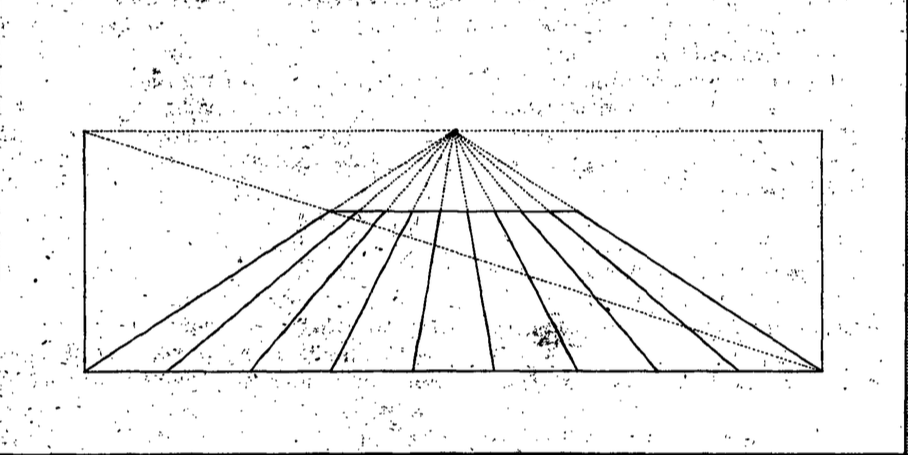}
{\centering\selectlanguage{english}
Figure 14. Leon Battista Alberti, \textit{Of Painting} \textit{in three books}, `Book I', in [2]\textit{.}
\par}

\bigskip
{\selectlanguage{english}
\textit{Step 2. Design the heights of the projections of the ``parallel'' straight lines of the floor} \textit{on the
painting }$\pi $\textit{. }\ The distance of the point of view of the painter from the vanishing point has to intervene
in this step. \ Consider the set painter-painting-floor seen from someone on the right, staying on the plane of the
painting $\pi $. Figure 15 shows how to construct these heights.}

\includegraphics[width=5.3736in,height=2.5827in]{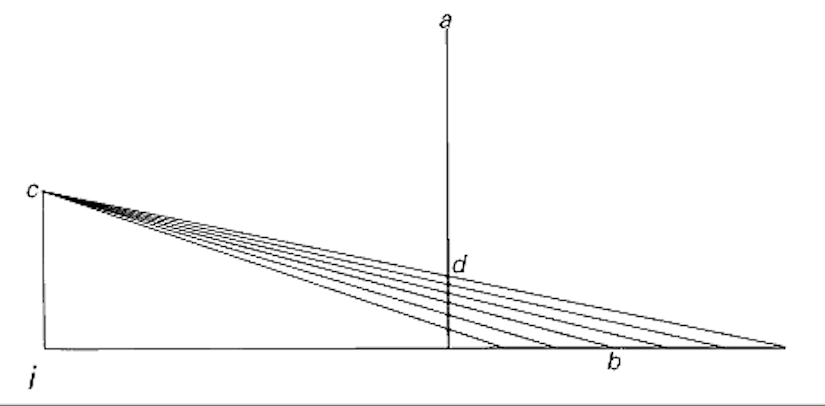}

{\centering\selectlanguage{english}
Figure 15
\par}

\bigskip
{\selectlanguage{english}
\textit{Step 3. Put together steps 1 and 2, and design the projection of the entire square-tile ground floor} \textit{on
the painting }$\pi $\textit{. } As shown in figure 16, it is enough to add to the painting obtained in Step 1 a
horizontal line at each of the heights constructed in Step 2.}

\includegraphics[width=5.6437in,height=2.3736in]{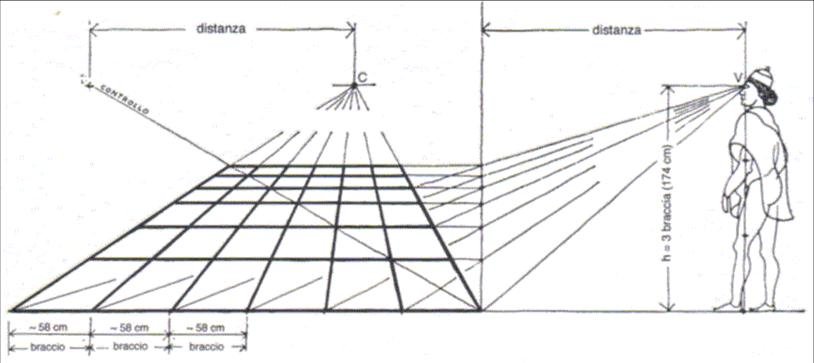}

{\centering\selectlanguage{english}
Figure 16
\par}

\bigskip
{\selectlanguage{english}
As the reader can see, the algorithmic construction illustrated by Leon Battista Alberti is very simple and does not
require any knowledge of sophisticated mathematical theories: this aspect made it really innovative at that time.}

{\selectlanguage{english}
Alberti, in his treatise \textit{De pictura praestantissima }[1]\textit{,} gives several other interesting and useful techniques for the
painters of his age, some of which are applications of the costruzione legittima. \ The process that we describe below
was meant to help the painter to identify the appropriate size of figures at different places on the square-tile floor,
exactly what Giotto would have needed in order to represent correctly the figure of San Francesco in the fresco we
described in Section 2. }

{\selectlanguage{english}
Note that the decision of placing the canvas on the floor, implies that only objects and figures placed on the floor
along the basis of the painting are represented in 1-1 scale. And one can then use the projection of the side of a
square tile parallel to the painting as a unit to give the measures of any object that is placed in the painting
precisely on this side (figure 17).}

\bigskip
{\centering 
\includegraphics[width=5.10in,height=3.20in]{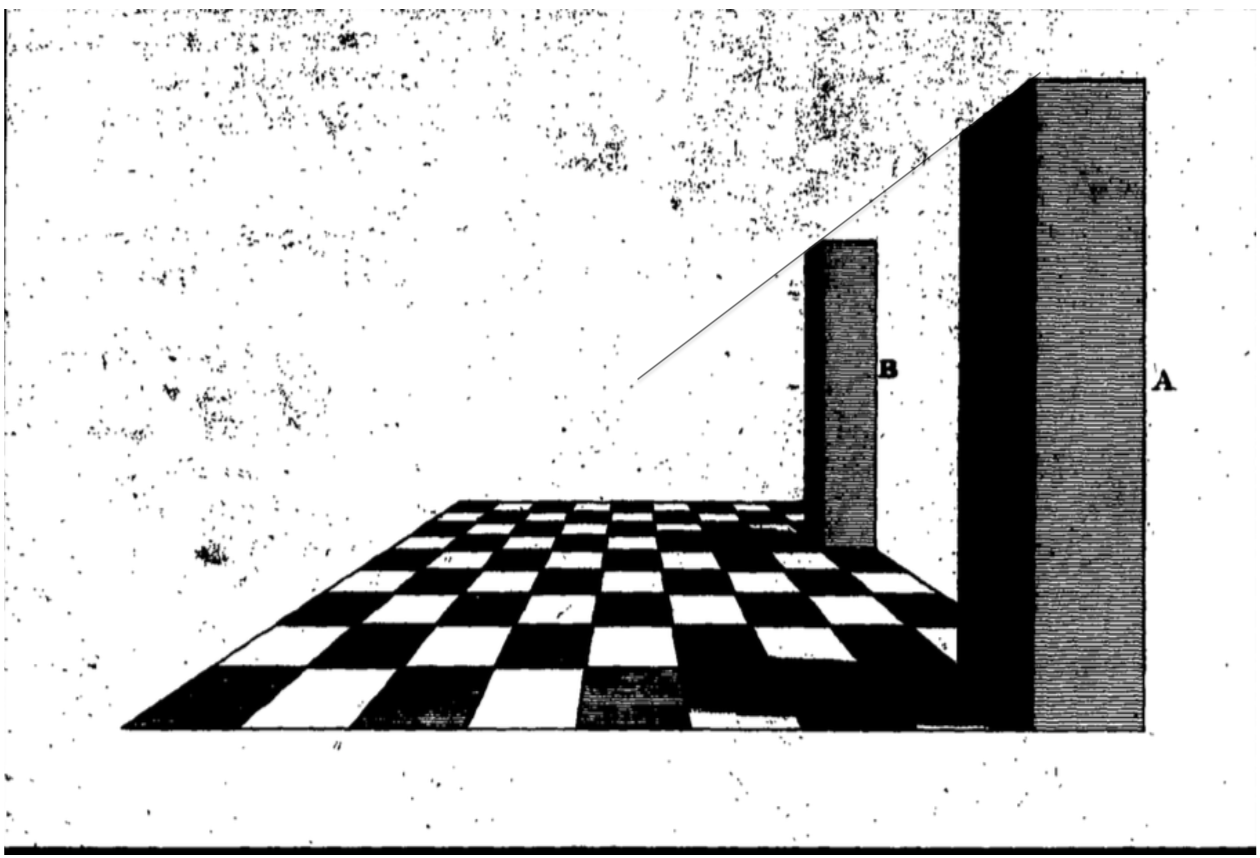}
\par}
{\centering\selectlanguage{english}
Figure 17. Leon Battista Alberti, \textit{Of Painting} \textit{in three books}, `Book II', in [2]\textit{.}
\par}

\bigskip
{\selectlanguage{english}
We now show how to use this representation to calculate the distance of the eye of the painter from the painting (figure
18). Note that a horizontal straight line L exiting from the eye and making an angle of 45 degrees with the plane of
the painting is parallel to one of the diagonals of the square tiles. Therefore, the line L, and all of the parallel
diagonals, encounter the horizon of the painting at a same point $A$ (see Section 3). \ Therefore, by extending
a diagonal of a tile in the painting until it encounters the horizon, one can find the point $A$. And now, since
the triangle with vertices the eye, the vanishing point, and the point $A$ is isosceles (it has the angles equal
to 45 degrees), then the distance between $A$ and the vanishing point is equal to the distance of the eye from
the painting. Therefore, by means of the side-of-tile-meter one can find the desired distance. \ }
\vskip .5cm 
{\centering 
\includegraphics[width=5.4957in,height=3.6264in]{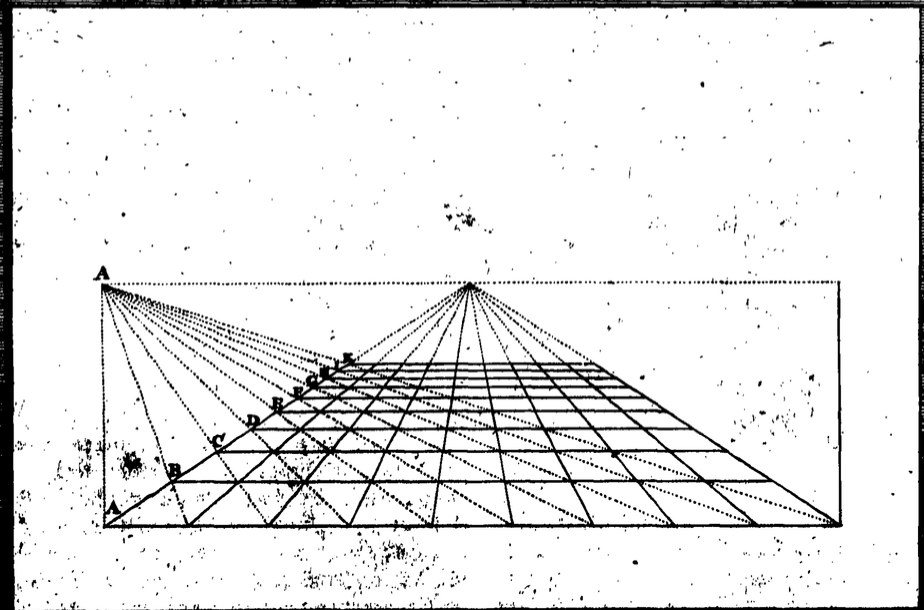}
\par}
{\centering\selectlanguage{english}
Figure 18. Leon Battista Alberti, \textit{Of Painting} \textit{in three books}, `Book I', in [2]\textit{.}
\par}
\bigskip
{\selectlanguage{english}
But Alberti, faithful to his promise to provide a practical manual and not a mathematical one, simply says to the
painter: \ extend one of the diagonals of a tile of the floor until it encounters the horizon of the painting at
$A$. Measure the distance between $A$ and the vanishing point. That distance is equal to the distance of
the eye of the painter from the painting itself. }

{\selectlanguage{english}
\section{5. Reconstructing a scene from a painting.}

{\selectlanguage{english}
As we have shown in the previous sections, the entire purpose of perspective is to take a three-dimensional scene, and
translate it into a two-dimensional scene (the painting) in a way that would fool the viewer into believing he is
actually looking at the original scene. }

{\selectlanguage{english}
But one could ask the inverse question. Can we reconstruct a real life scene, by just looking at its painting?
Immediately, we should know that the answer is negative. In fact it is clear that when we go from three dimensions down
to two dimensions we must lose some piece of information: the simplest way to convince ourselves of this consists in
closing one eye and start walking around. It will become easily apparent that a single eye provides a lack of depth
that may make some of the common chores difficult. This is essentially because the eye acts like a projection
mechanism, and the image of a three dimensional object is there represented as a flat picture on the bottom of the
retina. To remedy this difficulty, most animals have developed a system with two eyes.}

{\selectlanguage{english}
We could say that a properly designed painting is like a 2D compression of the data of a 3D scene. And, at least in
special situations, the originating scene can be appropriately reconstructed.}

{\selectlanguage{english}
The first situation in which a reconstruction is possible is the one in which a painting has a floor. If this is the
case, and if one knows both the distance $D$ of the point of view O from the plane of the painting ${\pi}$, and the
height H of the point of view from the floor, then all the vertical figures and objects that are standing on the floor
can be well placed in 3D. This is made clear by the following Thales-style\footnote{\textrm{ By this we mean a drawing
that utilizes a theorem that is often referred to as Thales' Theorem, namely an important result in elementary geometry
about the ratios of different line segments that arise if two intersecting lines are intercepted by two parallel
lines.}} drawings (figures 19, 20):}

\includegraphics[width=5.20in,height=2.75in]{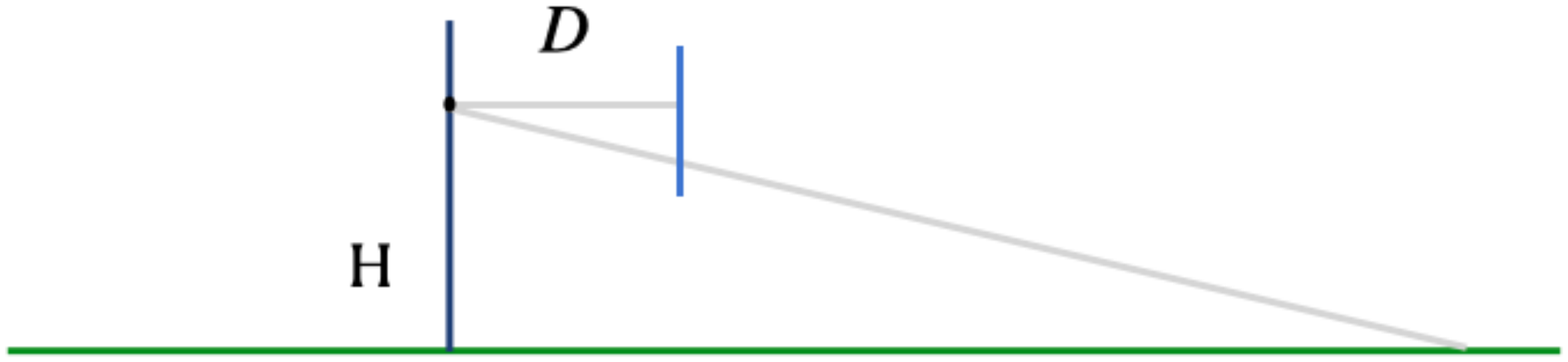}

{\centering\selectlanguage{english}
Figure 19
\par}

\bigskip

\includegraphics[width=5.20in,height=2.75in]{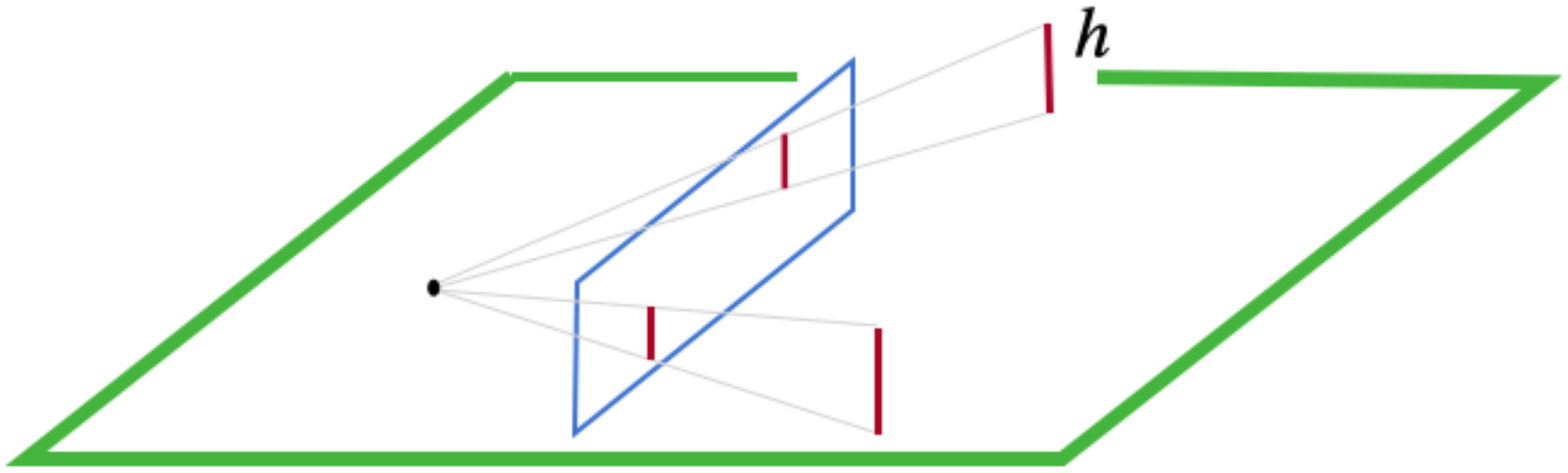}

{\centering\selectlanguage{english}
Figure 20
\par}

\bigskip
{\selectlanguage{english}
For those objects and figures that touch the floor, reconstruction is very easy: one has only to project from the point
of view through the painting and reach the ground, on the other side of the painting. Hence the height of each vertical
figure/object will become clear.}

{\selectlanguage{english}
But how can one recover the measures of $H$ and $D$, that are -- as it is clear -- fundamental for the reconstruction?
\ These key measures belong to the world that was external to the painting at the time it was painted: after several
centuries, the external world and the participating characters have all disappeared. The only hope is to find pieces of
information concerning that world, encoded inside the painting. It is like if we need to enter the painting, in a new
Mary Poppins-type walk. }

{\selectlanguage{english}
Let us now see how this was done in a specific case [16], for Piero della Francesca's \ \textit{Flagellazione }(figure
21), a first example -- we should say the example -- of the mathematically well constructed theory of perspective
contained in his \textit{De prospectiva pingendi}\footnote{ \textrm{In the extensive literature dedicated to the
}\textrm{\textit{Flagellazione }}\textrm{the article of Wittkower and Carter [27] offers an analysis made with
particular attention to the technical aspects of the perspective and to the historical language used. In this essay the
authors also trace the influence of the painter's perspective on real architecture.}\par }. }

\bigskip

{\centering 
\includegraphics[width=5.8835in,height=4.178in]{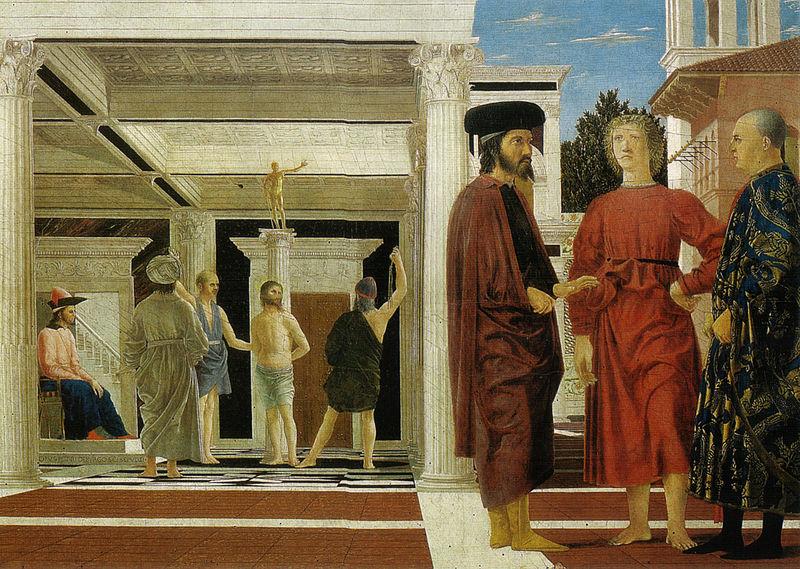}
\par}
{\centering\selectlanguage{english}
\foreignlanguage{italian}{Figure 21. Piero della Francesca, }\foreignlanguage{italian}{\textit{Flagellazione di Cristo,
}}\foreignlanguage{italian}{(1444-1470)}\foreignlanguage{italian}{\textit{, }}\foreignlanguage{italian}{tempera on
wood, Galleria Nazionale delle Marche, Urbino.}
\par}

\bigskip
{\selectlanguage{english}
First one notes that there are several human figures in the painting, whose knees all touch the line of the horizon
(figure 22); since, in the early renaissance and as we have mentioned before when discussing Alberti, the height of the
painted human figure was rigidly fixed to be three ``braccia fiorentine'', one immediately deduces that the height of
the knees, of the horizon, of the vanishing point $V$, and finally of point of view $O$ of the painter turns out to be
approximately 60 cm. }

{\centering  \includegraphics[width=5.6173in,height=4.061in]{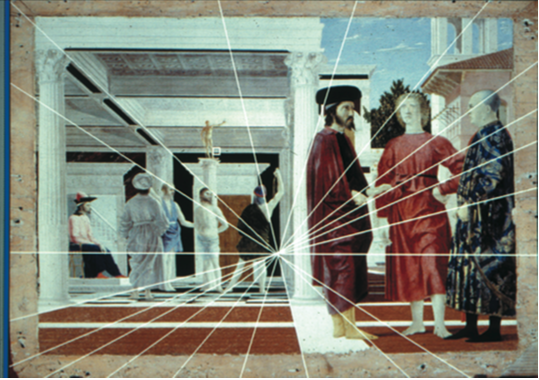}
\par}
{\centering\selectlanguage{english} Figure 22 \par}

\bigskip

{The determination of the distance $D$ is even more challenging. Here is how it was performed in the case of the \textit{Flagellazione}. 

As shown in figure 23, there is exactly one straight line L exiting from the point of view $O$ (the eye of the painter),
intersecting the horizon of the painting in a point $PD$ (the distance point, see Section 3) on the right side of the
vanishing point $PV$, and such that the triangle $O$, $PV$, $PD$, formed by the eye $O$, the vanishing point $PV$, and the point $PD$
is isosceles and rectangle in $PV$.

As we have seen in Section 3, all straight lines of the space that are parallel to
$L$, when represented in the painting, will have the same vanishing point $PD$. Therefore, if we find a line $K$ in the
painting ${\pi}$ that is the projection of a straight line of the space parallel to $L$, then we can solve the problem:
we intersect the extension of $K$ with the horizon and find the point $PD$, and then we try to figure out the ``real''
distance between $PV$ and $PD$, which will be the real distance between the eye $O$ of the painter and the painting ${\pi}$.}

{\selectlanguage{english}
Note that the floor of the \textit{Flagellazione} has a rectangular tile, whose real sides are (horizontal and)
parallel, respectively orthogonal, to the painting. If the tile of the floor were square, then one of the diagonals of
a tile would be a possible line $K$ we are searching for. 

{\selectlanguage{english}
Further more, we see a decoration of the floor, inscribed in a rectangular tile near the column with the Christ,
rendered as an ellipse in the painting. Of course this decoration could be in the reality either an ellipse or a
circle.\footnote{\textrm{ This is due to the fact the circle and the ellipse are two possible sections of the visual
cone that has vertex in the eye. In reality, and this goes beyond the purpose of this article, projective geometry
gives us all the tools to describe the way in which circles are transformed when we project them from the ceiling-floor
system to the painting. }} If it were a circle, then we could deduce that the tile is a square, and that its diagonal
is a possible line $K$. Then we could intersect its extension with the horizon and find $PD$, which will in turn suggest
the distance between $PV$ and $PD$, and hence the distance between the eye of the painter $O$ and the painting ${\pi}$. }

{\selectlanguage{english}
But now art history comes to our help. It appears that in the late 1400's no elliptical decorations were used in a floor
of tiles, and hence it can be now demonstrated that the distance between the point of view $O$ (the eye of the painter)
and the painting\textcolor{red}{ }is\textcolor{red}{ }approximately cm 145, [16]. }

{\centering \includegraphics[width=4.750in,height=4.00in]{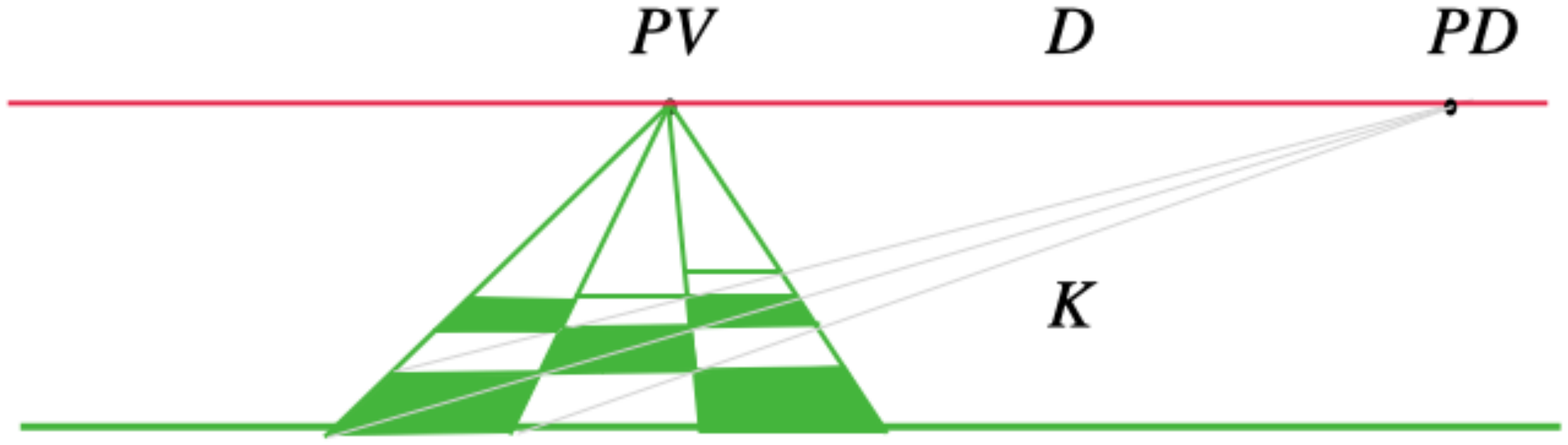}
\par}
{\centering \selectlanguage{english}
Figure 23
\par}
\bigskip

{\selectlanguage{english}
We conclude this section with a comment on the use of perspective not simply to represent reality, but to attribute
additional meanings to it. }

{\selectlanguage{english}
We believe that the mathematical analysis and reconstruction of the three dimensional scene of Piero's
\textit{Flagellazione} presented above could add a technical contribution to the historical-iconological one, in
connection with the hermeneutical problem that alimented the main interpretations of this painting proposed in the last
fifty years.\footnote{\textrm{ For a synthesis of the debate, see e.g., [13, p. 54 and ff].}} \ }

{\selectlanguage{english}
The identification of the figures of the painting, and in particular of the three of them which appear in the foreground
relies upon such a scant documentation, so that the iconological enigma hidden in the \textit{Flagellazione} seems to
remain still unsolved. In particular, the identification of the figure who appears on the right hand\footnote{\textrm{
As a curiosity we point out that the painting that is reproduced in [14] is actually a specular image of the actual
painting, a minor mistake that does not alter the interest of the article.}} side of the painting in the blue brocade
tunic -- likely an exponent of the noble Montefeltro family and possibly the patron of the painting -- remains
uncertain, [13, p. 62 and ff].}

{\selectlanguage{english}
The representation of patrons is not a surprising fact (akin to the naming of buildings that we see on campuses around
the world), but it is often somewhat unrelated to the painting itself. As an example, we can remind the reader of the
Scrovegni Chapel in Padua, where Giotto depicts the patron (Enrico Scrovegni) in the act of donating the chapel to the
Holy Virgin (figure 24). }

\bigskip

{\centering  \includegraphics[width=3.0083in,height=3.348in]{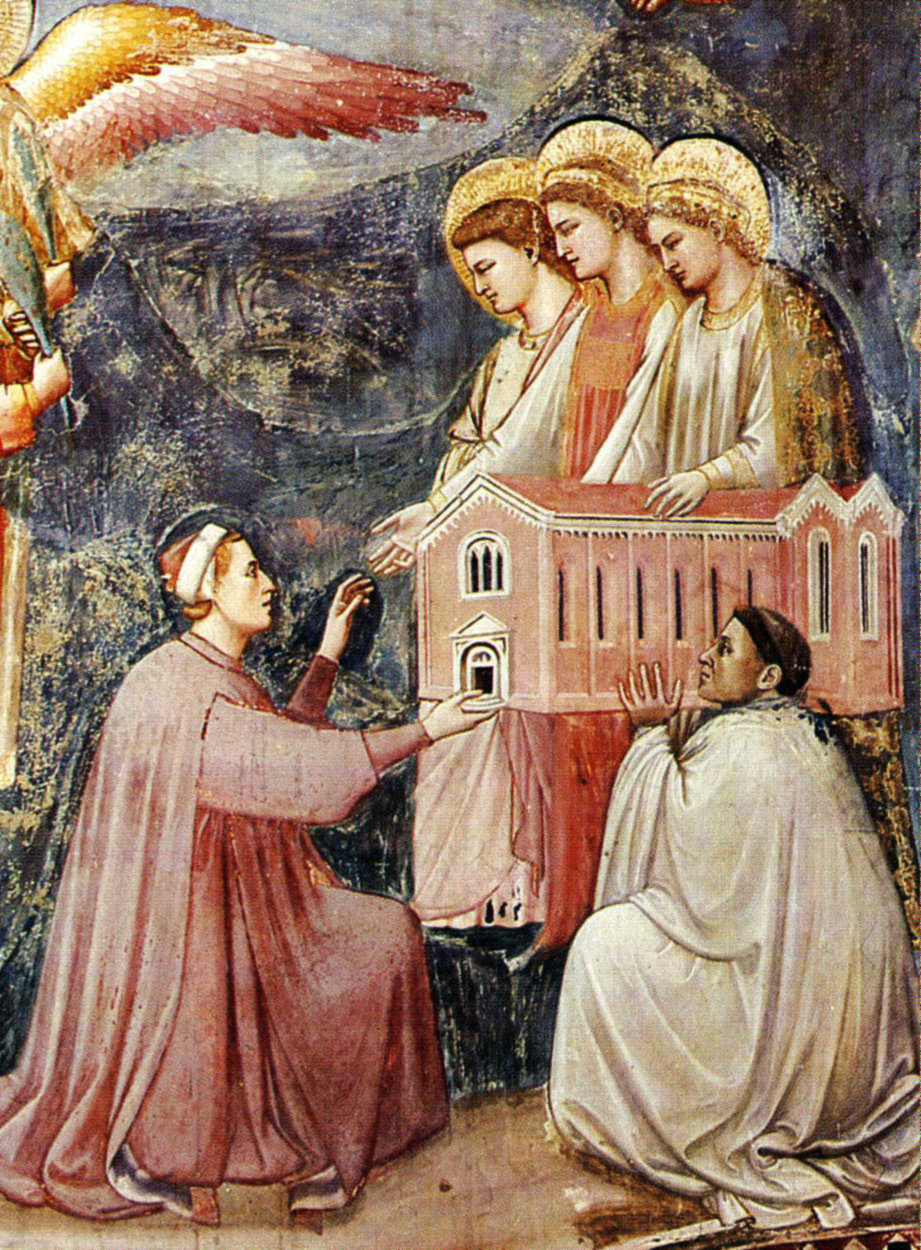}
\par}
{\centering \selectlanguage{english}
\foreignlanguage{italian}{Figure 24. Giotto, }\foreignlanguage{italian}{\textit{Last
Judgment}}\foreignlanguage{italian}{, ($\sim$1305), detail: Enrico Scrovegni gives to Madonna the model of Cappella
Scrovegni, affresco, Cappella degli Scrovegni, Padova. }
\par}

\bigskip
{\selectlanguage{english}
In Piero's case, however, we believe we can read an attempt to place the patron exactly at the scene, through the use of
perspective technique. Indeed, the fact that the level of the knees of the three gentlemen on the right is the same
level for the knees of Jesus and his torturers, indicates very specifically that the two sides of the picture were
rendered by the painter as if they were taking place at the same place and the same time. \ We see, therefore,
perspective used not simply as a geometrical device, but as a narrative instrument: the painter has made, here, a very
specific choice to insert the contemporary figures in a way that places them within the context of the historical
event.}

{\selectlanguage{english}
\section{6. Conclusions}

{\selectlanguage{english}
As we indicated in the introduction, this article is dedicated to the linear aspects of perspective. We have used the
desire of renaissance painters to faithfully represent tables, ceiling beams, and square floor decorations, to create
the new object that mathematician call the projective plane. This object (which will represent the floor and the
ceiling) is nothing but the old plane, to which one must add new (improper) points to represent the vanishing points
that the eye sees in a scene, as well as a new line, which is the line to which all improper points belong, and that is
represented as the horizon. But the story of perspective and projective geometry does not end here. The next natural
step, at least for a mathematician, is to study second degree equations, such as circles, ellipses, and other conic
sections. From the point of view of the painter, this is also an urgent matter, as it relates to the representation of
important everyday objects such as plates, carriage wheels, and windows, as well as not so everyday objects (yet very
important in religious paintings) such as halos. Projective geometry, with ideas that go back to the ancient Greek
mathematicians, provides a beautiful and very elegant solution to this problem, but this will be the object of a
subsequent article. }

%

%

\bigskip

\bigskip

{\selectlanguage{english}
\foreignlanguage{italian}{\textbf{\Large References}}}

{\selectlanguage{english}
\noindent \foreignlanguage{italian}{[1] Leon Battista Alberti, }\foreignlanguage{italian}{\textit{De pictura praestantissima, et
numquam satis laudata arte libri tres absolutissimi, Leonis Baptistae de Albertis viri in omni scientiarum genere, \&
praecipue mathematicarum disciplinarum doctissimi. }}\textit{Iam primum in lucem editi}, Westheimer, Basel, 1540.}

\bigskip

{\selectlanguage{english}
\noindent [2] Leon Battista Alberti, \textit{The Architecture {\dots} \ in ten books. Of Painting in three books. And Of Statuary
in one book. Translated into Italian by Cosimo Bartoli, and }\textit{into English by James Leoni, Architect.
Illustrated with seventy-five copper-plates, engraved by Mr. Picart,} Edward Owen, London, 1755.}

\bigskip

{\selectlanguage{english}
\noindent [3]\textit{ }Ibn al-Haytham,\textit{ The optics. Books 1-3 On direct vision, }translated with introduction and
commentary by A I Sabra, Warburg Institute, University of London, London, 1989.}

\bigskip

{\selectlanguage{english}
\noindent [4] Kirsti Andersen, \textit{The Geometry of an Art. The History of the Mathematical Theory of Perspective from Alberti
to Monge}, Springer, New York, 2007.}
\bigskip

{\selectlanguage{english}
\noindent [5] Hans Belting, \textit{Perspective: Arab Mathematics and Renaissance Western Art}, European Review 16, no. 2 (2008),
pp. 183-190.}

\bigskip

{\selectlanguage{english}
\noindent [6] Hans Belting, \textit{La double perspective. }\foreignlanguage{italian}{\textit{La science arabe et l'art de la
Renaissance}}\foreignlanguage{italian}{, La presse du reel/Presses universitaires de Lyon, Lyon, 2010.}}

\bigskip

{\selectlanguage{english}
\noindent [7] Hans Belting, \textit{The Double Perspective: Arab Mathematics and Renaissance Art}, Third Text 24, no. 5 (2010),
pp. 521-527.}

\bigskip

{\selectlanguage{english}
\noindent [8] Rudolf Bkouche, \textit{La naissance du projectif.} \textit{De la perspective à la géométrie projective}, in
\ Roshdi Rashed,~\textit{Mathématiques et Philosophie.~De l'Antiquité à l'Âge classique, }Paris,\textit{~}C.N.R.S.
Editions, 1991, pp. 239-285.}

\bigskip

{\selectlanguage{english}
\noindent [9] James Elkins, \textit{Piero della Francesca and the Renaissance. Proof of Linear Perspective,}
The Art Bulletin 69, no. 2 (1987), pp. 220-230.}

\bigskip

{\selectlanguage{english}
\noindent [10] J.V. Field, \textit{Alberti, the `Abacus' and Piero della Francesca's proof of perspective}, Renaissance Studies
11, no. 2 (1997), pp. 61-88.}

\bigskip

{\selectlanguage{english}
\noindent [11] J.V. Field, \foreignlanguage{english}{\textit{The Invention of Infinity. Mathematics and Art in the Renaissance}},
Oxford University Press, Oxford, New York, Tokyo, 1997.}

\bigskip

{\selectlanguage{english}
\noindent  \foreignlanguage{italian}{[12] J.V. Field, }\foreignlanguage{italian}{\textit{Piero della Francesca. }}\textit{A
Mathematician's Art}, Yale University Press, New Haven and London, 2005.}

\bigskip

{\selectlanguage{english}
\noindent  \foreignlanguage{italian}{[13] Carlo Ginzburg, }\foreignlanguage{italian}{\textit{Indagini su Piero. Il battesimo, il
ciclo di Arezzo, la Flagellazione di Urbino}}\foreignlanguage{italian}{. Nuova edizione, Einaudi, Torino, 1994.}}

\bigskip

{\selectlanguage{english}
\noindent  \foreignlanguage{italian}{[14] Martin Kemp, }\foreignlanguage{italian}{\textit{Piero's
perspective}}\foreignlanguage{italian}{, Nature}\foreignlanguage{italian}{\textit{ }}\foreignlanguage{italian}{390, no.
13 (1997), p. 128.}}

\bigskip

{\selectlanguage{english}
\noindent \foreignlanguage{italian}{[15] Piero della Francesca, }\foreignlanguage{italian}{\textit{De prospectiva
pingendi}}\foreignlanguage{italian}{, MSS }\url{http://digilib.netribe.it/bdr01/Sezione.jsp?idSezione=50}}

{\selectlanguage{english}
electronic reproduction: }

{\selectlanguage{english}
\noindent  \url{http://digilib.netribe.it/bdr01/visore/index.php?pidCollection=De-prospectiva-pingendi:889&v=-1&pidObject=De-prospectiva-pingendi:889&page=001%20R}}

\bigskip

{\selectlanguage{english}
\noindent  \foreignlanguage{italian}{[16] Placido Longo, }\foreignlanguage{italian}{\textit{La «Flagellazione» di Piero della
Francesca fra Talete e Gauss}}\foreignlanguage{italian}{, Bollettino dell'Unione Matematica Italiana 8, no. 2 (1999),
pp. 121--144. }\url{http://www.bdim.eu/item?id=BUMI_1999_8_2A_2_121_0}}

\bigskip

{\selectlanguage{english}
\noindent [17] Lucretius.~\textit{On the Nature of Things, }Translated by~W. H. D. Rouse.~Revised by~Martin F. Smith, Harvard
University Press, Cambridge, MA, 1924.}

\bigskip

{\selectlanguage{english}
\noindent [18] Erwin Panofsky, \textit{Perspective as symbolic form}, Zone Books, New York - MIT Press, Cambridge Mass., 1991.}

\bigskip

{\selectlanguage{english}
\noindent [19] Herman Schüling, \textit{Geschichte der Linear-Perspektive im Lichte der Forschung von ca 1870-1970},
Universitatsbibliothek, Giessen, 1975.}

\bigskip

{\selectlanguage{english}
\noindent [20] Gerard Simon, \textit{Optique et perspective: Ptolémée, Alhazen, Alberti / Optics and perspective: Ptolemy,
Alhazen, Alberti}, Revue d'histoire des sciences 54, no. 3 (2001), pp. 325-350.}

\bigskip

{\selectlanguage{english}
\noindent \foreignlanguage{italian}{[21] Luigi Vagnetti, `}\foreignlanguage{italian}{\textit{De naturali et artificiali
perspectiva':~bibliografia ragionata delle fonti teoriche e delle ricerche di storia della prospettiva. Contributo alla
formazione della conoscenza di un'idea razionale, nei suoi sviluppi da Euclide a Gaspard
Monge}}\foreignlanguage{italian}{, Edizione della Cattedra di composizione architettonica IA di Firenze e della L.E.F.,
1979.}}

\bigskip

{\selectlanguage{english}
\noindent \foreignlanguage{italian}{[22] Giorgio Vasari, }\foreignlanguage{italian}{\textit{Le vite de' più eccellenti pittori,
scultori et architettori}}\foreignlanguage{italian}{, Lorenzo Torrentino, Firenze, 1550, 3 voll.}}

{\selectlanguage{english}
\noindent \foreignlanguage{italian}{\ }\url{http://vasari.sns.it/vasari/consultazione/Vasari/ricerca.html}}

\bigskip

{\selectlanguage{english}
\noindent \foreignlanguage{italian}{[23] Giorgio Vasari, }\foreignlanguage{italian}{\textit{Live of Leon Battista
Alberti}}\foreignlanguage{italian}{ in: Vasari's }\foreignlanguage{italian}{\textit{Lives of the
Artists}}\foreignlanguage{italian}{. }\url{http://members.efn.org/~acd/vite/VasariAlberti.html}}
\bigskip

{\selectlanguage{english}
\noindent [24]\textit{ }Giorgio Vasari, \textit{Lives of the Artists}, translated with an Introduction and Notes by Julia Conaway
Bondanella and Peter Bondanella, Oxford University Press, Oxford, 1991.}

\bigskip

{\selectlanguage{english}\color{black}
\noindent [25] Kim H. Veltman, \textit{Literature on Perspective. A Select Bibliography (1971-1984)}, Marburger Jahrbuch für
Kunstwissenschaft 21, (1986), pp. 185-207.}

\bigskip

{\selectlanguage{english}
\noindent \foreignlanguage{italian}{[26] Graziella Federici Vescovini, }\foreignlanguage{italian}{\textit{De la métaphysique de la
lumière à la physique de la lumière dans la perspective des XIIIe et XIVe siècles}}\foreignlanguage{italian}{, Revue
d'histoire des sciences 60, no. 1 (2007), pp. 101-118.}}

\bigskip

{\selectlanguage{english}
\noindent [27] R. Wittkower and B. A. R. Carter, \textit{The Perspective of Piero della Francesca's `Flagellation'}, Journal of
the Warburg and Courtauld Institutes 16, no. 3/4 (1953), pp. 292-302.}
\bigskip

\begin{center}

\bigskip

{\large 
\noindent Dipartimento di Matematica e Informatica “U. Dini”, Universit\`a di Firenze,\\
Viale Morgagni 67/A, I-50134 Firenze, Italy.  graziano.gentili@unifi.it
\bigskip

\noindent Istituto per la storia del pensiero filosofico e scientifico moderno C.N.R.\\
 Area 3 - Bicocca Milano, via Cozzi, 53, 20125 Milano, Italy.  luisa.simonutti@ispf.cnr.it
 \bigskip
 
\noindent Donald Bren Presidential Chair in Mathematics, Chapman University,\\ 
One University Drive, Orange, CA 92866, USA.  struppa@chapman.edu
}
\end{center}

\end{document}